\documentclass[12pt]{article}

\usepackage[english]{babel}
\usepackage[utf8]{inputenc}

\usepackage{algorithm,color}
\usepackage[intlimits]{amsmath} 
\usepackage{amsthm,mathrsfs}    
\usepackage[comma,authoryear]{natbib}
\usepackage{url}
\numberwithin{equation}{section}
\usepackage{amssymb,amsmath}
\usepackage{subfigure}
\usepackage{bm}
\usepackage{booktabs} 
\usepackage{graphicx,graphics,epsfig}

\usepackage{float}
\usepackage{multirow}
\usepackage{array}

\theoremstyle{plain}
\newtheorem{thm}{Theorem}

\newcommand{\bthm}{\begin{thm}}
\newcommand{\ethm}{\end{thm}}

\newcommand{\bpf}{\begin{proof}}
\newcommand{\epf}{\end{proof}}

\theoremstyle{definition}
\newtheorem{defn}{Definition}
\newtheorem{rem}{Remark}
\newtheorem{motivation}{Motivation}
\newtheorem{note}{Note}
\newtheorem{sig}{Significance}

\DeclareMathOperator*{\argmin}{\arg\!\min}

\renewcommand{\baselinestretch}{1.2}
\usepackage[margin=.9in,a4paper]{geometry}
\usepackage[compact]{titlesec}
\titlespacing*{\section}{0pt}{1pt}{1pt}
\titlespacing*{\subsection}{0pt}{1pt}{1pt}
\titlespacing*{\subsubsection}{0pt}{1pt}{1pt}
\usepackage{deep-macro}
\begin{document}
\begin{center}
{\Large{\bf Unified Statistical Theory of Spectral Graph Analysis}} \\[.25in] 
Subhadeep Mukhopadhyay\\
Department of Statistical Science,  Temple University\\ Philadelphia, Pennsylvania, 19122, U.S.A.\\[1.25em]
\emph{Dedicated to the beloved memory of Emanuel (Manny) Parzen.}\\[1.6em]
\end{center}

\begin{abstract} 
The goal of this paper is to show that there exists a simple, yet universal\footnote{By universal, I mean the theory should not be specific to any particular spectral graph analysis technique but be meaningfully applicable to a class of models.} statistical logic of spectral graph analysis by recasting it into a nonparametric function estimation problem. The prescribed viewpoint appears to be good enough to accommodate most of the existing spectral graph techniques as a consequence of just one single formalism and algorithm.
\end{abstract} \vspace{-.4em}
\noindent\textsc{\textbf{Keywords and phrases}}: Nonparametric spectral graph analysis; Graph correlation density field (\texttt{GraField}); Empirical and smoothed spectral graph analysis; High-dimensional discrete data smoothing.
\vskip1.65em
\renewcommand{\baselinestretch}{.5}
\setlength{\parskip}{.25ex}
{\small
\tableofcontents
}
\renewcommand{\baselinestretch}{1.2}
\setlength{\parskip}{1.25ex}
\section{Introduction} 
Spectral graph data analysis is undoubtedly the most favored technique for graph data analysis, both in theory and practice. It has been extraordinarily successful in many applications including data mining, web search, quantum computing, computer vision, image segmentation, and others. Over the last several decades, extensive efforts have been made by researchers and practitioners to develop a suit of spectral graph analysis techniques (e.g., Laplacian, Modularity, Diffusion map, regularized Laplacian, Google PageRank model etc.) with increasing sophistication and specialization. However, no single algorithm can be regarded as a panacea for dealing with the evolving complexities of modern graphs. Therefore, the most important and pressing question for the field today appears to be whether we can develop a unifying language to establish ``bridges'' between a wide variety of spectral graph analysis techniques, and thus providing logically connected means for reaching different ends. Undoubtedly, any such formulation would be of great theoretical significance and practical value, that will ultimately provide the applied data scientists clear guidance and a systematic strategy for selecting the proper spectral tools to arrive at a confident conclusion.
\vskip.66em

To that end, this work attempts to unify the theories of spectral graph analysis by purely statistical means. That is, we seek a constructive (rather than confirmatory) theoretical framework integrating classical and modern spectral graph analysis algorithms, which to date have been viewed as distinct methods.
\subsection{Spectral Graph Analysis: A Practitioner’s Guide} \label{sec:prac}
The way spectral graph analysis is currently taught and practiced can be summarized as follows (also known as \emph{spectral heuristics}):
\vspace{-.4em}
\begin{enumerate}
  \item Let $\cG=(V,E)$ denotes a  (possibly weighted) undirected graph with a finite set of vertices $|V|=n$, and a set of edges $E$. Represent the graph using weighted adjacency matrix $A$ where $A(x,y;\cG)=w_{xy}$ if the nodes $x$ and $y$ are connected by an edge and $0$ otherwise; weights are non-negative and symmetric.
  \item Define ``suitable'' spectral graph matrix (also known as graph ``shift'' operator). Most popular and successful ones are listed below:
  \begin{itemize}
    \item $\cL=D^{-1/2}AD^{-1/2}$;  ~~~\cite{chung1997}
    \item $\mathcal{B}=A \,-\,N^{-1} d d^{T}$; ~~~\,\cite{newman2006}
    \item $\calT=D^{-1}A$; ~~~\,\cite{coifman06}
    \item Type-I Reg. $\cL_\tau= $ $D_\tau^{-1/2}A\,D_\tau^{-1/2}$;\, ~~~\,\cite{chaudhuri2012}
    \item Type-II Reg. $\cL_\tau=D_\tau^{-1/2}\,A_\tau\,D_\tau^{-1/2}$; ~~\cite{amini2013}
    \item Google's PageRank $\calT_\al=\al D^{-1}A + (1-\al) F$; ~~\citep{GOOGLE}
  \end{itemize} \vspace{-.2em}
  $D = {\rm diag}(d_1,\ldots,d_n) \in \cR^{n \times n}$, $d_i$ denotes the degree of a node, $\tau>0$ regularization parameter, $F\in \cR^{n \times n}$ with all entries $1/n$, and $N=2|E|=\sum_{x,y}A(x,y)$.
  \item Perform spectral decomposition of the matrix selected at step 2, whose eigenvectors form an orthogonal basis of $\cR^n$. Spectral graph theory seeks to understand the interesting properties and structure of a graph by using leading nontrivial eigenvectors and eigenvalues, first recognized by \cite{fiedler73}, which provide the basis for performing harmonic analysis on graphs.
   \item Compute graph Fourier transform (GFT) by expanding functions (or signals) defined over the vertices of a graph as a linear combination of the selected eigenbasis (from step 3) and carry out learning tasks such as regression, clustering, classification, smoothing, kriging, etc.
\end{enumerate}
\vspace{-.6em}
\subsection{Previous Theoretical Treatments}
There has been a flood of publications on statistical spectral graph theory,  most of which have focused on asymptotic justifications of different techniques (mentioned in the previous section). See, in particular, those based on Laplacian \citep{von2008consistency,belkin2001laplacian}, modularity \citep{davis2016}, diffusion map \citep{coifman06, coifman2005pnas}, and regularized Laplacian  \citep{chaudhuri2012,le2015sparse} methods. More refined theoretical results in the context of a specific parametric graph model (such as stochastic block model or its variants) are discussed in \cite{joseph2013, rohe2011spectral,sarkar2015} and \cite{qin2013}. 
\vskip.55em

\emph{What is this paper not about}? It is important to draw a distinction between the theory that purports to provide the constructions of different spectral graph techniques and, the theory that does not. This article focuses on the former. Instead of statistical \emph{confirmation}, here our main interest lies in understanding the  statistical \emph{origins} of various spectral graph techniques. This is mostly ill-understood and neglected territory. The current theoretical treatments do not give us any clues on this issue.

\subsection{A Set of Questions} 
Despite many advances, the theory and practice of spectral graph analysis are still very much compartmentalized and unsystematic. The graph signal processing engineers are often left bewildered by the vastness of the existing literature and huge diversity of methods (developed by the machine learning community, applied harmonic analysts, physicists, and statisticians). This has led to a need to develop a broader perspective on this topic, lest we be guilty of not seeing the forest for the trees. The question \emph{how} different spectral graph techniques can naturally originate from some underlying basic principle, plays the key role in this regard by clarifying the mystery of \emph{why and when} to use them.  The reality is we still lack general tools and techniques to attack this question in a statistical way. Given the very different character of the existing spectral models, it seems that new kinds of abstractions need to be formulated to address this fundamental challenge. What would these abstractions look like?  How can we discover a \emph{statistical path} that naturally leads to these different spectral methods? One such promising theory is discussed in this paper, taking inspiration from celebrated history on nonparametric spectral analysis of time series. The prescribed modern viewpoint shows an exciting confluence of nonparametric function estimation, quantile domain methods, and Fourier harmonic analysis to unveil a more simple conceptual structure of this subject, avoiding undue complexities.

\subsection{Unified Construction Principle}
We start by giving a glance to the unified spectral representation scheme, whose (nonparametric) statistical interpretation and justification is deferred until the next section, as it requires a more technical background. It turns out that all known spectral graph techniques are just different manifestations of this single general algorithm in some way.
\vskip.5em
\emph{Probability Notations and Definitions.} We will stick to the following probabilistic notations. Define network probability mass function by $ P(x,y;\cG)=A(x,y;\cG)/\sum_{x,y}A(x,y;\cG)$; Vertex probability mass function by $p(x;\cG)=\sum_{y} P(x,y;\cG)$ and $p(y;\cG)=\sum_{x} P(x,y;\cG)$ with the associated quantile functions $Q(u;X,\cG)$ and $Q(v;Y,\cG)$. Finally, define the important Graph Interaction Function by ${\rm GIF}(x,y;\cG)\,=\,P(x,y;\cG)/p(x;\cG) p(y;\cG),$ $x,y\in V(\cG)$.

\begin{center}
\emph{Algorithm 1. Nonparametric Spectral Graph Analysis: A Unified Algorithm}
\end{center}
\vspace{-.4em}
\medskip\hrule height .65pt
\vskip.8em

\texttt{Step 1.} For given discrete graph $\cG$ of size $n$, construct GraField kernel function $\cd: [0,1]^2 \rightarrow \cR_+ \cup \{0\}$ defined a.e by
\beq \label{eq:gcddef}
\cd(u,v;\cG)\,=\,{\rm GIF}\big[Q(u;X),Q(v;Y);\cG\big]\,=\,\dfrac{P\big( Q(u;X), Q(v;Y);\cG \big)}{p\big( Q(u;X) \big) p\big(  Q(v;Y) \big)},
\eeq
where $u=F(x;\cG),v=F(y;\cG)$ for $x,y \in \{1,2,\ldots,n\}$.
\vskip1em
\texttt{Step 2.}  Let $\{\xi_j\}_{j\ge0}$ be a complete orthonormal system on $[0,1]$. Construct discrete transform basis $\xi_j(x;F(X;\cG)):=\xi_j(F(x;\cG))$, evaluated at the vertex-\emph{rank-transforms}. They satisfy the following orthonormality conditions (degree-weighted):
\[\sum_x \xi_j(x;F(X)) \,\xi_k(x;F(X))\, p(x;\cG) \,= \,\delta_{jk},\]
thus constitutes an orthonormal basis for $L^2(F;\cG)$.
\vskip1em
\texttt{Step 3.} Transform coding of graphs. Construct generalized spectral graph matrix $\mathcal{M}(\cG,\xi)\in \cR^{n\times n}$ with respect to $\eta_j(u):=\xi_j(Q(u;X);F)$:
\beq \label{eq:ogt} \mathcal{M}[j,k;\cG, \xi]\,=\Big\langle \eta_j, \int_0^1(\cd-1)\eta_k \Big\rangle_{L^2[0,1]}=\,\sum_{\ell,m} \xi_j(\ell;F) \xi_k(m;F) P(\ell,m;\cG).\eeq
$\mathcal{M}_\xi$ can be viewed as a \emph{transform coefficient matrix} of the orthogonal series expansion of $\cd(u,v;\cG)$ with respect to the product bases $\{\eta_j\eta_k\}_{1 \le j,k \le n}$.
\vskip1em
\texttt{Step 4.} Perform the singular value decomposition (SVD) of $\mathcal{M}_\xi=U\Lambda U^{T}$ $= \sum_k u_k \mu_k u_k^{T}$, where where $u_{ij}$ are the elements of the singular vector of moment matrix $U=(u_1,\ldots,u_n)$, and $\Lambda={\rm diag}(\mu_1,\ldots,\mu_n)$, $\mu_1 \ge$ $ \cdots \mu_n \ge 0$.
\vskip1em
\texttt{Step 5.} Obtain approximate \KL (KL) representation basis (which acts as a graph Fourier basis) of the graph $\cG$ by
\[\widetilde{\phi}_{k}\,=\,\sum_{j=1}^n u_{jk}\xi_{j},\,~ {\rm for}\,\, k=1,\ldots,n-1,\]
which can be directly used for subsequent signal processing on graphs.
\vskip.55em
\medskip\hrule height .65pt
\vskip1.5em
\emph{Organization.} The next section is devoted to the fundamental principles that underlie this algorithm. The unifying power will be shown in Sections 3 and 4. By doing so, we will introduce empirical and smooth spectral graph analysis techniques. We also address the open problem of obtaining a formal interpretation of ``spectral regularization''. A deep connection between high-dimensional discrete data smoothing and spectral regularization is discovered. This new perspective provides, for the first time, the theoretical motivation and fundamental justification for using regularized spectral methods, which were previously considered to be empirical guesswork-based ad hoc solutions. Important application towards spatial graph regression is discussed in Section 5. In Section 6, we end with concluding remarks. Additional applications are available as supplementary materials.
\section{Fundamentals of Statistical Spectral Graph Analysis}
\subsection{Graph Correlation Density Field}
Wiener's generalized harmonic analysis formulation on spectral representation theory of time series \emph{starts} by defining the autocorrelation function (ACF) of a signal. In particular,  Wiener–Khintchine theorem asserts ACF and the spectral density are Fourier duals of each other.  Analogous to Wiener's \emph{correlation} technique, we describe a new promising starting point, from which we can develop the whole spectral graph theory systematically by bringing nonparametric function estimation and harmonic analysis perspectives. Both statistical and probabilistic motivations will be given.


\begin{defn}
For given discrete graph $\cG$ of size $n$, the piecewise-constant bivariate kernel function $\cd: [0,1]^2 \rightarrow \cR_+ \cup \{0\}$
is defined almost everywhere through
\beq \label{eq:gcddef}
\cd(u,v;\cG_n)\,=\,{\rm GIF}\big[Q(u;X),Q(v;Y);\cG_n\big]\,=\,\dfrac{p\big( Q(u;X), Q(v;Y);\cG_n \big)}{p\big( Q(u;X) \big) p\big(  Q(v;Y) \big)},~~~0<u,v<1.
\eeq
\end{defn}
\begin{thm} \label{lemma1} GraField defined in Eq. (\ref{eq:gcddef}) is a positive piecewise-constant kernel satisfying
\[\iint_{[0,1]^2} \cd(u,v;\cG_n) \dd u \dd v\,~ = \sum_{(i,j)\in \{1,\ldots,n\}^2} \iint_{I_{ij}} \cd(u,v;\cG_n) \dd u \dd v ~= ~1,\]
where
\[I_{ij}(u,v)= \left\{ \begin{array}{ll}
         \,1, ~~& {\rm if}~ (u,v) \in \left(F(i;X), F(i+1;X)\right] \times \left( F(j;Y), F(j+1;Y)\right]\\
         \,0, ~~& \mbox{elsewhere}.\end{array} \right.\]
\end{thm}
\vskip.25em
\begin{note}
The bivariate step-like shape of the \texttt{GraField} kernel is governed by the (piecewise-constant left continuous) quantile functions $Q(u;X,\cG_n)$ and $Q(v;Y,\cG_n)$ of the \emph{discrete} vertex probability measures. As a result, in the continuum limit (let $(\cG_n)_{n\ge0}$ be a sequence of graphs whose number of vertices tends to infinity $n \rightarrow \infty$), the shape of the piecewise-constant discrete $\cd_n$ approaches a ``continuous field'' over unit interval -- a self-adjoint compact operator on square integrable functions (defined on the graph) with respect to vertex probability measure $p(x;\cG)$ (see Section 2.2 for more details).
\end{note}
\vskip.15em
\begin{motivation}
We start with a diffusion based probabilistic interpretation of the GraField kernel. The crucial point to note is that the ``slices'' of the $\cd$ \eqref{eq:gcddef} can be expressed as $p(y|x;\cG)/p(y;\cG)$ in the vertex domain. This alternative conditional probability-based viewpoint suggests a connection with the random walk on the graph. Interpret $p(y|x;\cG)$ as transition probability from vertex $x$ to vertex $y$ in one time step. Also note that $p(y;\cG)$ is the stationary probability distribution on the graph, as we have $\lim_{t\rightarrow \infty} p(t,y|x;\cG)=p(y;\cG)$ regardless of the initial starting point $x$ (moreover, for connected graphs the stationary distribution is unique). Here $p(t,y|x;\cG)$ denotes the probability distribution of a random walk landing at location $y$ at time $t$, starting at the vertex $x$. See \cite{lovasz1993random} for an excellent survey on the theory of random walks on graphs.

The graph affinity function measures how the transition probability $p(y|x;\cG)$ is different from the ``baseline'' stationary distribution (long-run stable behavior) $p(y;\cG)$. That \emph{comparison ratio} is the fundamental interaction function for graphs which we denote by ${\rm GIF}(x,y;\cG)$. This probabilistic interpretation along with Theorem 2 and 5 will allow us to integrate the \emph{diffusion map} \citep{coifman06} technique into our general statistical framework in Section 3.2.
\end{motivation}
\vskip.1em
\begin{motivation}
GraField compactly represents the affinity or strength of ties (or interactions) between every pair of vertices in the graph.
To make this clear, let us consider the following adjacency matrix of a social network representing $4$ employees of an organization
\[ A=\left( {\large\begin{smallmatrix}
0 ~~ &2 ~ ~ &0  ~~ &0\\[.25em]
2 ~ ~&0  ~~ &3  ~~ &3\\[.25em]
0 ~~ &3  ~~ &0  ~~ &3\\[.25em]
0 ~ ~&3  ~~ &3 ~ ~ &0
\end{smallmatrix}}\right), \]
where the weights reflect numbers of communication (say email messages or coappearances in social events etc.). Our interest lies in understanding the strength of association between the employees i.e., ${\rm Strength}(x,y)$ for all pairs of vertices. Looking at the matrix $A$ (or equivalently based on the histogram network estimator $p(x,y;\cG)=A/N$ with $N=\sum_{x,y} A(x,y)=22$) one might be tempted to conclude that the link between employee $1$ and $2$ is the weakest, as they have communicated only twice, whereas employees $2,3,$ and $4$ constitute strong ties, as they have interacted more frequently. Now, the surprising fact is that (i) ${\rm Strength}(1,2)$ is twice that of  ${\rm Strength}(2,3)$ and  ${\rm Strength}(2,4)$; also (ii) ${\rm Strength}(1,2)$ is $1.5$ times of ${\rm Strength}(3,4)$! To understand this paradox, compute the vertex-domain empirical GraField kernel matrix (Definition 1) with $(x,y)$th entry $N {\cdot} A(x,y;\cG)/d(x)d(y)$
\[ \cd_n=\left( {\large\begin{smallmatrix}
0  ~&22/8 ~  &0  ~ &0\\[.25em]
22/8  ~&0 ~  &22/16  ~ &22/16\\[.25em]
0 ~ &22/16  ~ &0  ~ &22/12\\[.25em]
0 ~ &22/16 ~ &22/12 ~  &0
\end{smallmatrix}}\right). \]
This toy example is in fact a small portion (with members $1,9,31$ and $33$) of the famous Zachary's karate club data, where the first two members were from  Mr. Hi’s group and the remaining two were from John’s group. The purpose of this illustrative example is not to completely dismiss the adjacency or empirical graphon \citep{lovasz2006} based analysis but to caution the practitioners so as not to confuse the terminology ``strength of association'' with ``weights'' of the adjacency matrix -- the two are very different objects. Existing literature uses them interchangeably without paying much attention. The crux of the matter is: association does \emph{not} depend on the raw edge-density, it is a ``comparison edge-density'' that is captured by the GraField; see Section 3.2 for its intriguing connection with diffusion based graph-distance.
\end{motivation}
\vskip.1em
\begin{motivation}
GraField can also be viewed as properly ``renormalized Graphon,'' which is reminiscent of Wassily Hoeffding's ``standardized distributions'' idea \citep{Hoeff40}. Thus, it can be interpreted as a discrete analogue of copula (the Latin word copula means ``a link, tie, bond'') density for random graphs that captures the underlying correlation field. We study the structure of graphs in the spectral domain via this fundamental graph kernel $\cd$ that characterizes the implicit \emph{connectedness} or \emph{tie-strength} between pairs of vertices.
\end{motivation}
Fourier-type spectral expansion results of the density matrix $\cd$ are discussed in the ensuing section, which is at the heart of our approach. We will demonstrate that this correlation density operator-based formalism provides a useful perspective for spectral analysis of graphs that allows unification.

\subsection{\KL Representation of Graph}
We define the \KL (KL) representation of a graph $\cG$ based on the spectral expansion of its \texttt{GraField} function $\cd(u,v;\cG)$. Schmidt decomposition \citep{schmidt1907} of $\cd$ yields the following spectral representation theorem of the graph.
\begin{thm} The square integrable graph correlation density kernel $\cd:[0,1]^2 \rightarrow \cR_+ \cup \{0\}$ of two-variables admits the following canonical representation
\beq
\cd(u,v;\cG_n)\,=\, 1 \,+\, \sum_{k=1}^{n-1} \la_k \phi_k(u) \phi_k(v),
\eeq
where the non-negative $\la_1 \ge \la_2\ge \cdots \la_{n-1}\ge 0$ are singular values and $\{\phi_k\}_{k\ge 1}$ are the orthonormal singular functions $\langle \phi_j,\phi_k\rangle_{\sL^2[0,1]}=\delta_{jk},$ for $j,k =1,\ldots,n-1$, which can be evaluated as the solution of the following integral equation relation
\beq \label{eq:IE}
\int_{[0,1]} \left[\cd(u,v;\cG) -1\right] \phi_k(v) \dd v\,=\,\la_k \phi_k(u),~~~~k=1,2,\ldots,n-1.
\eeq
\end{thm}

\begin{rem} By virtue of the properties of \KL(KL) expansion \citep{loeve1955book}, the eigenfunction basis $\phi_k$ satisfying \eqref{eq:IE} provides the optimal low-rank representation of a graph in the mean square error sense. In other words, $\{\phi_k\}$ bases capture the graph topology in the \emph{smallest} embedding dimension and thus carries practical significance for graph compression. Hence, we can call those functions the \emph{optimal} coordinate functions or Fourier representation bases. Accordingly, the fundamental statistical modeling problem hinges on finding approximate solutions to the optimal graph coordinate system $\{\phi_1,\ldots,\phi_{n-1}\}$ satisfying the integral equation \eqref{eq:IE}.
\end{rem}

\begin{defn} Any function or signal $y \in \cR^{n}$  defined on the vertices of the graph $y:V \mapsto \cR$ such that $\| y\|^2=\sum_{x \in V(\cG)} |y(x)|^2 p(x;\cG) < \infty$, can be represented as a linear combination of the Schmidt bases of the GraField density matrix $\cd$. Define the \emph{generalized graph Fourier transform} of $y$
\[\widehat{y}(\la_k)~:=~ \langle y, \phi_k \rangle ~=~ \sum_{x=1}^n y(x) \phi_k[F(x;\cG)].\]
This spectral or frequency domain representation of a signal, belonging to the square integrable Hilbert space $\sL^2(F,\cG)$ equipped with the inner product \[\large\langle y,z \large\rangle_{\sL^2(F;\cG)} = \sum_{x \in V(\cG)} y(x) z(x) p(x;\cG),\] allows us to construct efficient graph learning algorithms. As $\{\phi_k\}$'s are KL spectral bases, the vector of projections onto this basis function decay rapidly, and hence may be truncated  aggressively to capture the structure in a small number of bits.
\end{defn}

\begin{defn} The entropy (or energy) of a discrete graph $\cG$, is defined using the Parseval relation of the canonical representation
\[{\rm Entropy}(\cG)~=~\iint_{[0,1]^2} (\cd-1)^2\dd u \dd v ~= ~\sum_k |\la_k|^2.\]
This quantity, which captures the departure of uniformity of the $\cd$, can be interpreted as a measure of `structure' or the `compressibility' of the graph. Entropy measure can also be used to (i) define graph homogeneity and (ii) design fast algorithms for graph isomorphism. For homogeneous graphs, the shape of the correlation density field is uniform over the unit square. The power at each harmonic components, as a function of frequency, is called the power spectrum of the graph. Flat spectrum is equivalent to the analogous notion of white-noise process.
\end{defn}
\vspace{-.5em}
\subsection{Nonparametric Spectral Approximation Theory}
We view the \emph{spectral graph learning algorithm} as a method of approximating  $(\la_k,\phi_k)_{k\ge1}$ that satisfies the integral equation \eqref{eq:IE}, corresponding to the graph kernel $\cd(u,v;\cG)$. In practice, often the most important features of a graph can be well characterized and approximated by a few top (dominating) singular-pairs. The statistical estimation problem can be summarized as follows:
\[\vspace{-.4em}~~~~A_{n \times n} ~~\mapsto~~\cd ~~\mapsto~~ \left\{ \big(\widehat\la_1,\widehat\phi_1\big),\ldots,  \big(\widehat\la_{n-1},\widehat\phi_{n-1}\big) \right\}~\text{that satisfies Eq.}~\eqref{eq:IE}.~~~~~~\]

\vspace{-.4em}
\subsubsection{From Continuous to Discrete Basis and Back}
The Fourier-type nonparametric spectral approximation method starts by choosing an expansion basis. Let $\xi_j$ denote an orthonormal basis of $\sL^2[0,1]$. Construct the discrete transform basis of $\cR^n$ by evaluating $\xi_j$ at the vertex-rank-transforms $\xi_j(F(x;\cG)):=\xi_j(x;F,\cG)$. Verify that $\xi_j(x;F,\cG)$ are orthogonal with respect to the measure $p(x;\cG)$ satisfying
\[\sum_x \xi_j(x;F,\cG) \,\xi_k(x;F,\cG)\, p(x;\cG)~=~0,~~ \text{for}\, j \ne k.\]
Thus form an orthogonal basis of the Hilbert space $\sL^2(F;\cG)$. Define $\eta$-functions (quantile domain unit bases), generated from mother $\xi_j$ by
\[\eta_j(u;\cG)~=~\xi_j\big[ Q(u;X);F,\cG\big],~~~0<u<1,~~\]
as piecewise-constant (left-continuous) functions over the irregular grid $\{0,p(1), p(1)+p(2),\ldots, $ $\sum_{j=1}^n p(j)=1\}$ satisfying $\langle \eta_j, \eta_k\rangle_{\sL^2[0,1]}=0$, if $j\ne k$. They will provide a useful tool to recast conventional matrix calculus-based approaches as a functional statistical problem. One unique aspect of our construction is that in the continuum limit (as the size of the graph $n \rightarrow \infty$) the discrete $\eta$-basis of $\cR^n$ approaches the mother $\xi$-function, a basis of $\sL^2[0,1]$. We call it asymptotic ``reproducing'' property.
\subsubsection{Projection Methods for Eigenvector Approximation}
We are interested in the nonparametric estimation of eigenpairs $\{\la_k,\phi_k\}_{k\ge 0}$. Approximate the unknown eigenvectors by the projection, $\mathcal{P}_n\phi_k$, on the  \texttt{span}$\{\eta_j,j=1,\ldots,n\}$ defined by
\beq \label{eq:lc} \phi_k(u) \,\approx\, \mathcal{P}_n\phi_k\,=\, \sum_{j=1}^n \te_{jk}\, \eta_j(u),~~~0<u<1 \eeq
where $\te_{jk}$ are the unknown coefficients to be estimated.

\begin{defn}[Orthogonal Discrete Graph Transform] We introduce a generalized concept of matrices associated with graphs called the \texttt{G-matrix}. Define discrete graph transform with respect to an orthonormal system $\eta$ as
\beq \label{eq:ogt} \mathcal{M}[j,k; \eta, \cG] ~= ~\Big\langle \eta_j, \int_0^1(\cd-1)\eta_k \Big\rangle_{\sL^2[0,1]}~~~ \text{for}~ j,k=1,\ldots,n.\eeq
Equivalently, we can define the discrete graph transform to be the coefficient matrix of the orthogonal series expansion of the \texttt{GraField} kernel $\cd(u,v;\cG)$ with respect to the product bases $\{\eta_j\eta_k\}_{1 \le j,k \le n}$. As a practical significance, this generalization provides a systematic recipe for converting the graph problem into a ``suitable'' matrix problem:
\[~~~\cG_n(V,E)~~~ \longrightarrow~~~ A_{n\times n}~~~~ \longrightarrow ~ ~~~\cd(u,v;\cG_n) ~~~~ \xrightarrow[\,\,{\rm Eq.}\, (2.5)]{\{\eta_1,\ldots,\eta_n\}} ~ ~~~ \mathcal{M}(\eta,\cG)\,\in\, \cR^{n\times n}.~~~~~~\]
\end{defn}
\vspace{-.8em}
\begin{thm}
The \texttt{G-matrix} \eqref{eq:ogt} can also be interpreted as a ``covariance'' operator for a discrete graph by recognizing the following equivalent representation for $x,y\in \{1,2,.\ldots,n\}$
\[\mathcal{M}[j,k; \xi,\cG] ~= ~\Ex_P\left[ \xi_j\big( F(X;\cG) \big)  \xi_k\big( F(Y;\cG) \big) \right]~=~\sum_{x,y} P(x,y;\cG)\ \xi_j\big( F(x;\cG) \big) \ \xi_j\big( F(y;\cG) \big),\]
\end{thm} \vspace{-.25em}
This can be proved using the basic quantile mechanics fact that $Q(F(X)) = X$ holds with probability $1$ (see \cite{parzen79}). Next, we present a general approximation scheme that provides an effective method of discrete graph analysis in the frequency domain. \vskip1em
\begin{thm}[Nonparametric spectral approximation] \label{thm:GEN}
The Fourier coefficients $\{\te_{jk}\}$  of the projection estimators \eqref{eq:lc} of the GraField eigenfunctions (eigenvalues and eigenvectors), satisfying the integral equation \eqref{eq:IE},  can be obtained by solving the following generalized matrix eigenvalue problem
\beq \label{eq:master} \mathcal{M} \Te\,=\,S \Te \Delta,\eeq
where $\mathcal{M}_{jk}=\big\langle \eta_j, \int_0^1(\cd_n-1)\eta_k \big\rangle_{\sL^2[0,1]}$, and $S_{jk}=\big\langle \eta_j, \eta_j \big\rangle_{\sL^2[0,1]}$.
\end{thm}
To prove define the residual of the governing equation \eqref{eq:IE} by expanding $\phi_k$ as series expansion \eqref{eq:lc},
\beq \label{eq:er}
R(u) ~\equiv~ \sum_j \te_{jk} \Big[ \int_0^1 \big(  \cd(u,v;\cG) -1 \big) \eta_j(v) \dd v \,-\, \la_k  \eta_j(u)   \Big]\, =\, 0.
\eeq
Now for complete and orthonormal $\{\eta_j\}$ requiring the error $R(u)$ to be zero is equivalent to the statement that $R(u)$ is orthogonal to each of the basis functions
\beq \big\langle R(u),\eta_k(u)\big\rangle_{\sL^2[0,1]}~=~0,~~~~~~k=1,\ldots,n.  \eeq
This leads to the following set of equations:
\beq \label{eq:m1}
\sum_j \te_{jk} \Big[ \iint\limits_{[0,1]^2} \big(  \cd(u,v;\cG_n) -1 \big) \eta_j(v) \eta_k(u) \dd v \dd u\Big]~-~\la_k \sum_j \te_{jk} \big[\int_0^1 \eta_j(u) \eta_k(u) \dd u \big]~=~0.
\eeq
Theorem \ref{thm:GEN} plays a key role in our statistical reformulation. In particular, we will show how the fundamental equation \eqref{eq:master} provides the desired unity among different spectral graph techniques by systematically constructing a ``suitable'' class of coordinate functions.
\vskip.1em
\begin{note} The fundamental idea behind the Rietz-Galerkin \citep{galerkin1915} style approximation scheme for solving variational problems in Hilbert space played a pivotal inspiring role to formalize the statistical basis of the proposed computational approach.
\end{note}
\begin{note}
Our nonparametric spectral approximation theory based on eigenanalysis of \texttt{G-Matrix}, remains \emph{unchanged} for any choice $\eta$-function (following the recipe for construction given in Sec 2.3.1), which paves the way for the generalized harmonic analysis of graphs.
\end{note}
The next two sections investigate this general scheme under various choices of $\eta$-functions, and nonparametric estimation methods of $\cd_n$. By doing so, many ``mysterious similarities'' among different spectral graph algorithms are discovered which were not known before.

\section{Empirical Spectral Graph Analysis} 
Three popular traditional spectral graph analysis models will be synthesized in this section.
\subsection{Laplacian Spectral Analysis}
Laplacian is probably the most heavily used spectral graph technique in practice. Here we will demonstrate for the first time how the Laplacian of a graph \emph{naturally originates} by purely statistical reasoning, totally free from the classical combinatorial based logic.

\vskip.4em

{\bf Degree-Adaptive Block-pulse Basis Functions}.  One of the fundamental, yet universally valid (for any graph) choice for $\{\eta_j\}_{1\le j \le n}$ is the indicator top hat functions (also known as  block-pulse basis functions, or in short BPFs). However, instead of defining the BPFs on a uniform grid (which is the usual practice) here (following Sec 2.3.1) we define them on the non-uniform mesh $0=u_0 < u_1\cdots< u_n=1$ over [0,1], where $u_j=\sum_{x \le j}p(x;X)$ with local support
\beq \label{eq:basis}
\eta_j(u) = \left\{ \begin{array}{ll}
         p^{-1/2}(j) ~~& \mbox{for $u_{j-1} < u \le u_{j}$};\\
         0 ~~& \mbox{elsewhere}.\end{array} \right.
\eeq
They are disjoint, orthogonal, and a complete set of functions satisfying
\[\int_0^1 \eta_j(u) \dd u =\sqrt{p(j)},~~\,\int_0^1\eta_j^2(u) \dd u =1,~~\text{and}~~\int_0^1\eta_j(u)\eta_k(u) \dd u = \delta_{jk}.\]

\begin{figure*}[tth]
\centering
\includegraphics[height=\textheight,width=.41\textwidth,keepaspectratio,trim=3cm 1.5cm 2cm 3cm]{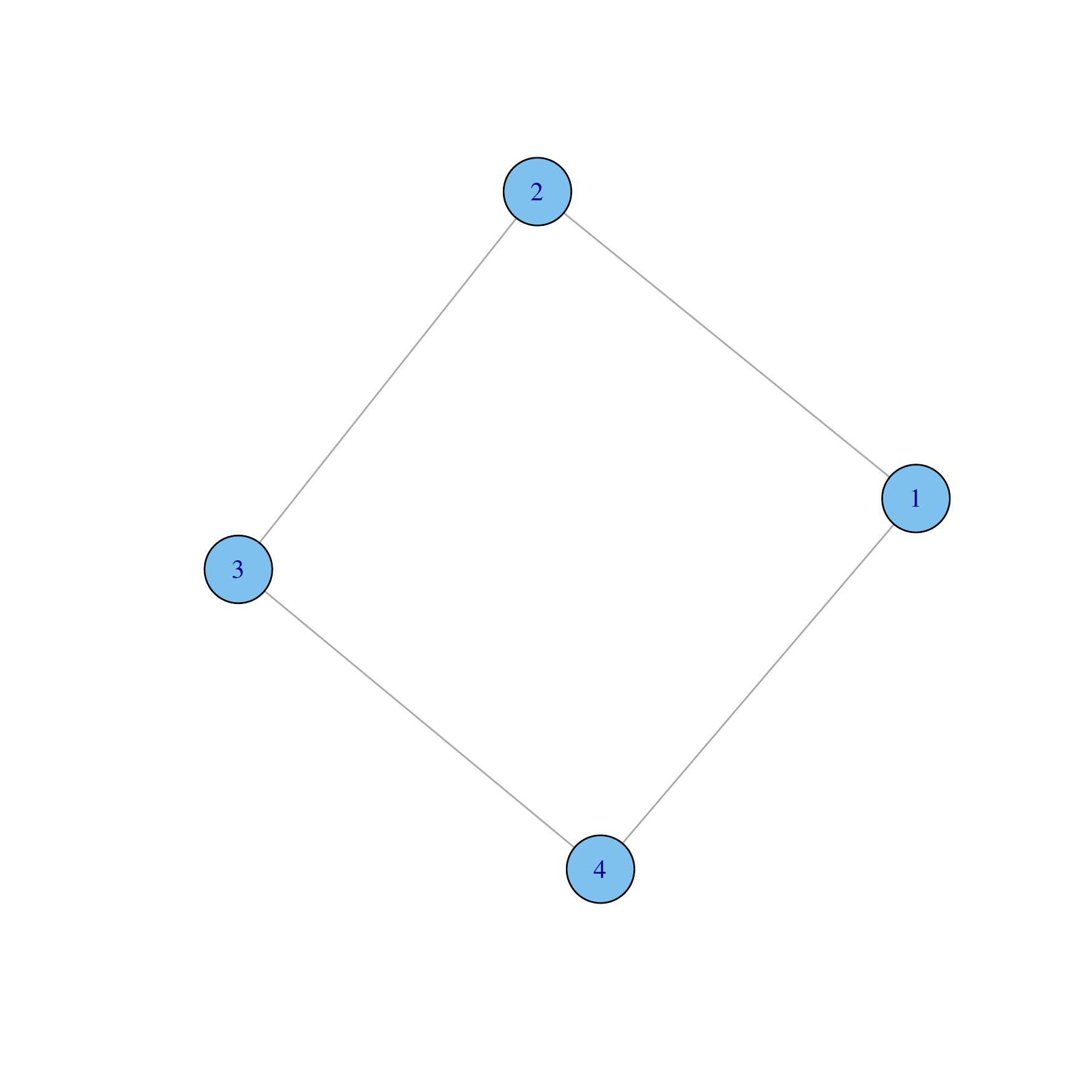}~~~~
\includegraphics[height=\textheight,width=.41\textwidth,keepaspectratio,trim=2cm 1.5cm 3cm 3cm]{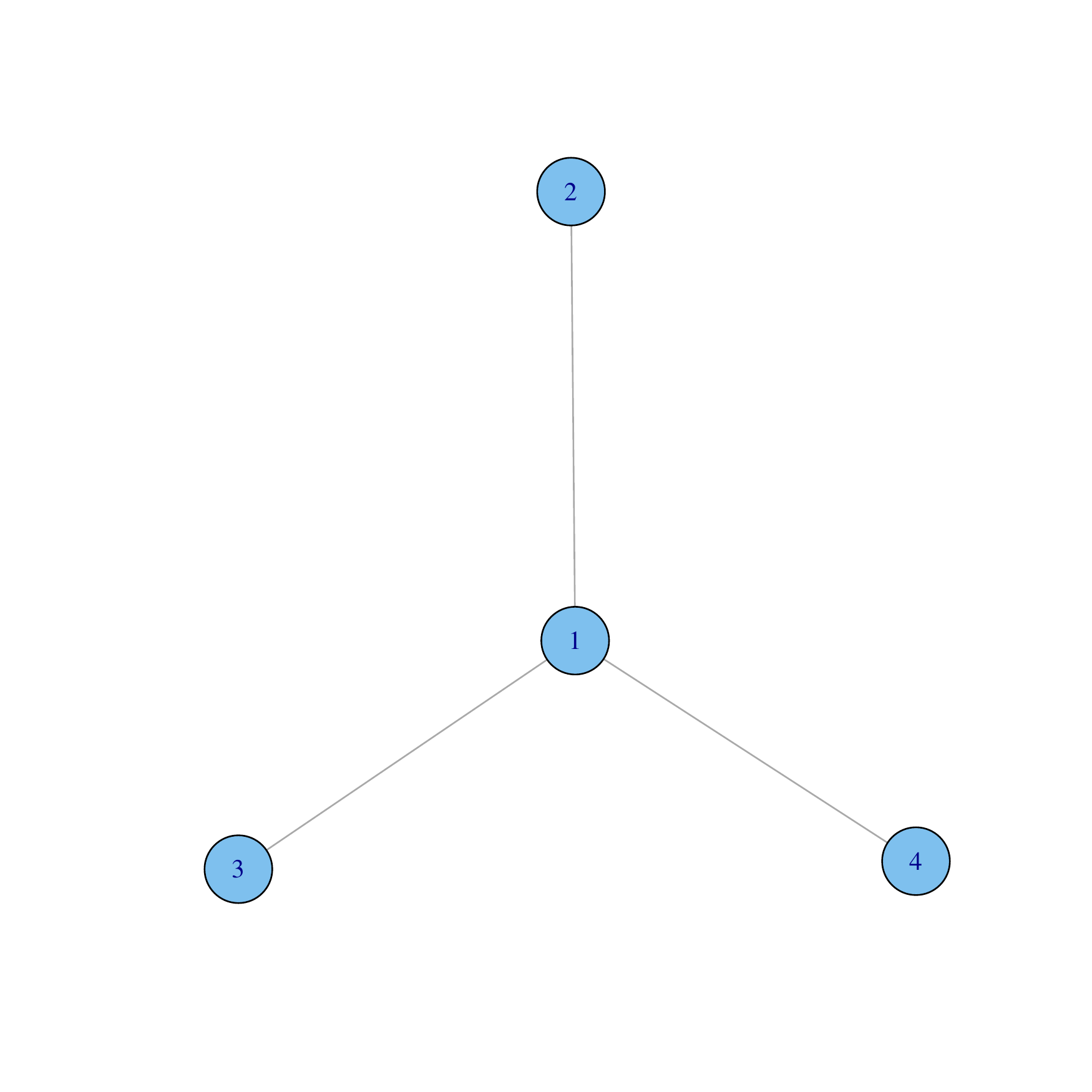}\\[1.5em]
\includegraphics[height=\textheight,width=.41\textwidth,keepaspectratio,trim=3cm 1cm .5cm 2cm]{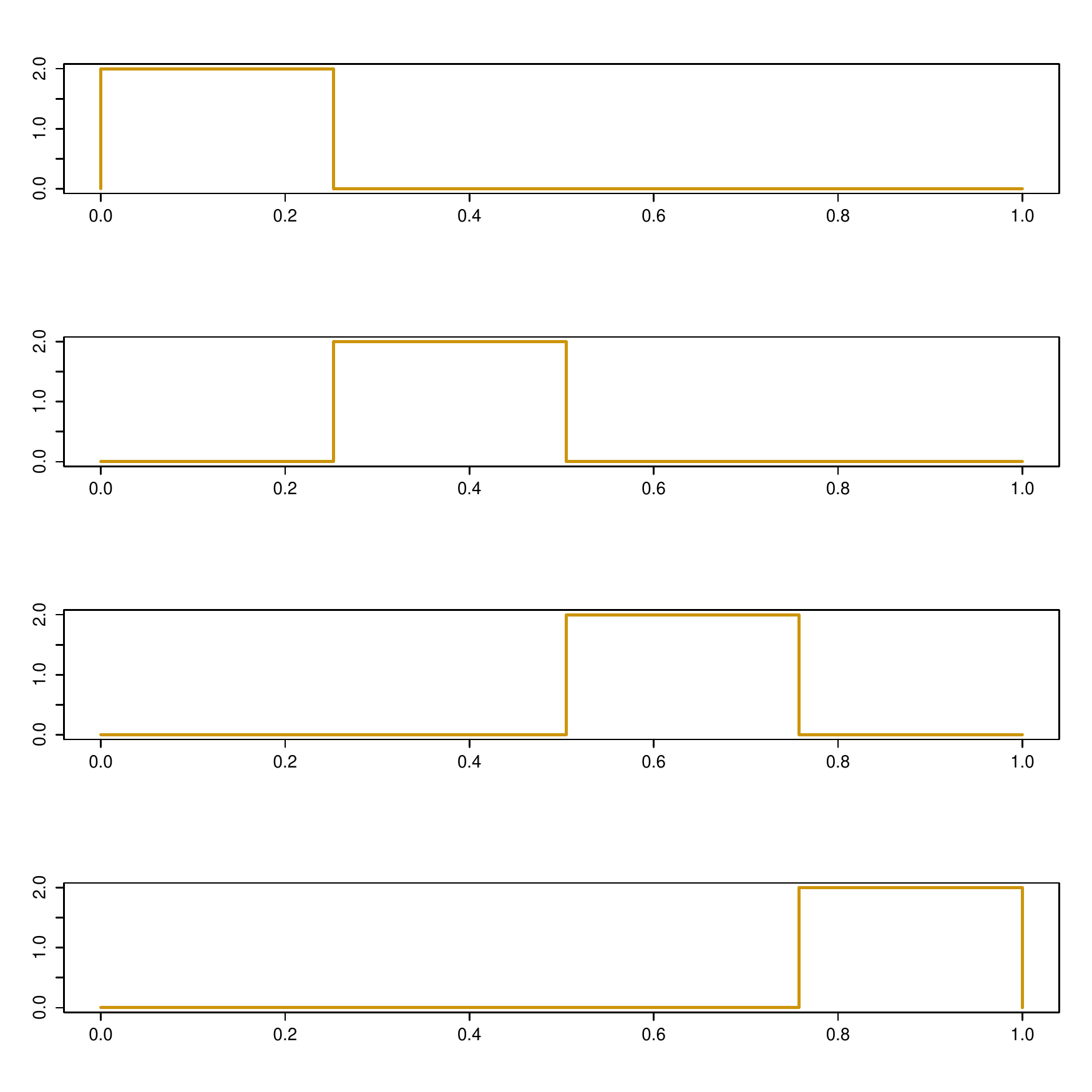}~~~
\includegraphics[height=\textheight,width=.41\textwidth,keepaspectratio,trim=.5cm 1cm 3cm 2cm]{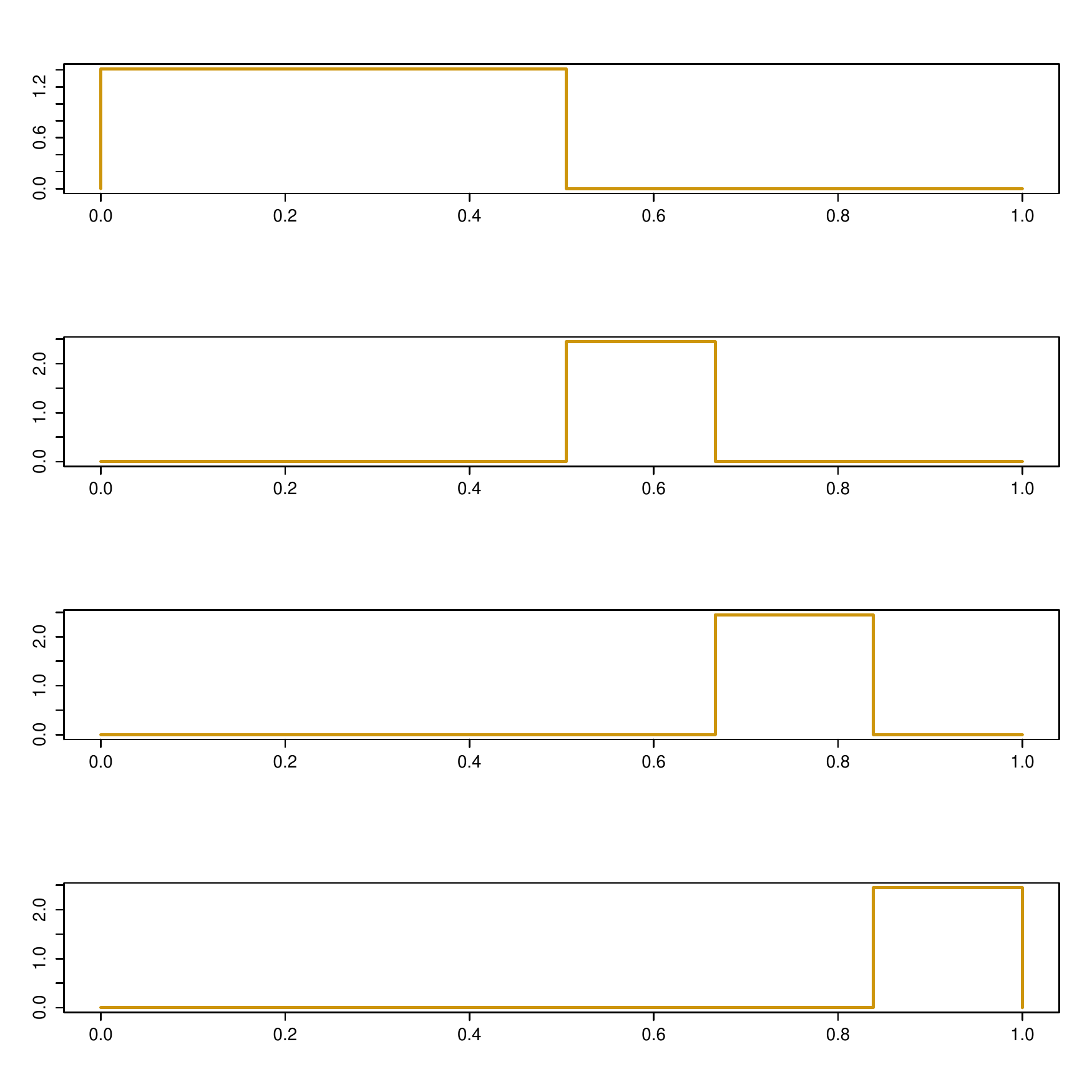}\\[1em]
\caption{Two graphs and their corresponding degree-adaptive block-pulse functions. The amplitudes and the block length of the indicator basis functions \eqref{eq:basis} depends on the degree distribution of that graph.} \label{fig:bpf}
\vspace{-.6em}
\end{figure*}

\begin{note} The shape (amplitudes and block lengths) of our specially designed BPFs depend on the specific graph structure via $p(x;\cG)$ as shown in Fig \ref{fig:bpf}. In order to obtain the spectral domain representation of the graph, it is required to estimate the spectra of GraField kernel $\phi_k$, by representing them as block pulse series. The next result describes the required computational scheme for estimating the unknown expansion coefficients $\{\te_{jk}\}$.
\end{note}

\begin{thm} \label{thm:lap}
Let $\phi_1,\ldots,\phi_n$ the canonical Schmidt bases of $\sL^2$ graph kernel $\cd(u,v;\cG)$, satisfying the integral equation \eqref{eq:IE}.
Then the empirical solution of \eqref{eq:IE} for block-pulse orthogonal series approximated \eqref{eq:basis} Fourier coefficients $\{\te_{jk}\}$ can equivalently be written down in closed form as the following matrix eigen-value problem
\beq \label{eq:ev}
\cL^*[\te]\,=\,\la \te,
\eeq
where $\cL^*=\cL-uu^{T}$, $\cL$ is the Laplacian matrix, $u=D_p^{1/2}1_n$, and $D_p={\rm diag}(p_1,\ldots,p_n)$.
\end{thm}

Note that the discrete GraField kernel takes value ${\rm GIF}(j,k;\cG_n)=P(j,k;\cG_n)/p(j;\cG_n)p(k;\cG_n)$ and the tensor–product bases $\eta_j(u)\eta_k(v)$ take value $p^{-1/2}(j;\cG_n)p^{-1/2}(k;\cG_n)$ over the rectangle $I_{jk}(u,v)$ for $0<u,v<1$. This observation reduces the master equation \eqref{eq:m1} to the following system of linear algebraic equations expressed in the vertex domain:
\beq \label{eq:mm2}
\sum_j \te_{jk} \left[ \dfrac{p(j,k)}{\sqrt{p(j)p(k)}}\,-\,\sqrt{p(j)}\sqrt{p(k)}\,-\,\la_k\delta_{jk} \right] ~=~0.
\eeq
\vskip.25em
{\bf Empirical plugin nonparametric estimator}. The estimating equation \eqref{eq:mm2} contains unknown network and vertex probability mass functions which need to be estimated from the data. The most basic nonparametric estimates are $\widetilde{P}(j,k;\cG_n)=A[j,k;\cG_n]/N$ and $\widetilde{p}(j;\cG_n)=d(j;\cG_n)/N$.  By plugging these empirical estimators into Eq. \eqref{eq:mm2}, said equation can be rewritten as the following compact matrix form:
\beq \label{eq:m3}
\left[\cL-N^{-1}\sqrt{d} \sqrt{d}^{T}\right] \tilde{\te} = \tilde{\la} \tilde{\te},
\eeq
where $\cL=D^{-1/2}AD^{-1/2}$ is the graph Laplacian matrix.
\vskip1em

\begin{sig} Theorem \ref{thm:lap} allows us to interpret the graph Laplacian as the empirical \texttt{G-matrix} $\widetilde{\mathcal{M}}(\eta,\cG_n)$ under degree-adaptive indicator basis choice for the $\eta$ shape function. Our technique provides a completely nonparametric statistical derivation of an algorithmic spectral graph analysis tool. The plot of $\widetilde{\la}_k$ versus $k$ can be considered as a ``raw periodogram'' analogue for graph data analysis.
\end{sig}

\subsection{Diffusion Map}
We provide a statistical derivation of Coifman's diffusion map \citep{coifman06} algorithm, which hinges upon the following key result.
\begin{thm} \label{thm:dm}
The \emph{empirical GraField} admits the following vertex-domain spectral \emph{diffusion} decomposition at any finite time $t$
\beq \label{eq:RW} \dfrac{\widetilde p(t,y|x;\cG_n)}{\widetilde p(y;\cG_n)}~=~1+\sum_k \tilde \la_k^t \tilde{\phi}_k(x;F) \tilde{\phi}_k(y;F),\eeq
where $\widetilde{\phi}_k=  D_{\tilde p}^{-1/2}u_k$, ($u_k$ is the $k$th eigenvector of the Laplacian matrix $\cL$), $(\phi_k \circ F)(\cdot)$ is abbreviated as $\phi_k(\cdot;F)$, $\widetilde p(y|x;\cG)=\calT(x,y)$, and $\calT=$ $D^{-1}A$ is the transition matrix of a random walk on $\cG$ with stationary distribution $\tilde{p}(y;\cG_n)=d(y;\cG_n)/N$.
\end{thm}
Replacing the estimated coefficients from Theorem \ref{thm:lap} into \eqref{eq:lc} yields $\widetilde{\phi}_k=  D_p^{-1/2}u_k$, where $u_k$ is the $k$th eigenvector of the Laplacian matrix $\cL$, immediately leads to the following vertex-domain spectral decomposition result of the empirical GraField: \beq \label{eq:RW} \dfrac{\widetilde p(y|x;\cG)}{\widetilde p(y;\cG)}~=~1+\sum_k \tilde \la_k \tilde{\phi}_k(x) \tilde{\phi}_k(y),\eeq
For an alternative proof of the expansion \eqref{eq:RW} see Appendix section of \cite{coifman06} by noting $\widetilde{\phi}_k$ are the right eigenvectors of random walk Laplacian $\calT$.
\begin{sig}
In light of Theorem \ref{thm:dm}, define diffusion map coordinates at time $t$ as the mapping from $x$ to the vector
\[ x \longmapsto \Big( \la_1^t\phi_1(x;F), \ldots,  \la_k^t\phi_k(x;F)\Big),~~~~x\in\{1,2,\ldots,n\},\]
which can be viewed as an approximate `optimal' \KL representation basis and thus can be used for non-linear embedding of graphs. Define diffusion distance, a measure of similarity between two nodes of a graph, as the Euclidean distance in the diffusion map space
\beq \label{eq:DD} D_t^2(x,x')~=~\sum_{j\ge 1} \la_j^{2t}\big\{\phi_j(x;F)  \,-\,\phi_j(x';F)\big\}^2.~~~~~~~~~~~\eeq
This procedure is known as \emph{diffusion map}, which has been extremely successful tool for manifold learning\footnote{Manifold learning: Data-driven learning of the ``appropriate'' coordinates to identify the intrinsic nonlinear structure of high-dimensional data. We claim the concept of GraField allows decoupling of the geometrical aspect from the probability distribution on the manifold.}. Our approach provides an \emph{additional} insight and justification for the diffusion coordinates by interpreting it as the strength of connectivity profile for each vertex, thus establishing a close connection with empirical GraField.
\end{sig}
\subsection{Modularity Spectral Analysis}
\begin{thm}
To approximate the KL graph basis $\phi_k= \sum_j\te_{jk} \eta_j$, choose $\eta_j(u)= \I(u_{j-1} < u \le u_{j})$ to be the characteristic function satisfying
\[\int_0^1 \eta_j(u) \dd u ~=~ \int_0^1 \eta^2_j(u) \dd u ~=~\,p(j;\cG). \]
Then the corresponding empirically estimated spectral graph equation \eqref{eq:m1} can equivalently be reduced to the following  generalized eigenvalue equation in terms of the matrix $\mathcal{B}=A \,-\,N^{-1} d d^{T}$
\beq
\mathcal{B}\al\,=\,\la D \al.
\eeq
\end{thm}

\begin{sig}  The matrix $\mathcal{B}$, known as modularity matrix, was introduced by \cite{newman2006} from an entirely different motivation. Our analysis reveals that the Laplacian and Modularity based spectral graph analyses are equivalent in the sense that they inherently use the same underlying basis expansion (one is a rescaled version of the other) to approximate the optimal graph bases.
\end{sig}

\vskip.4em

\begin{rem} Solutions (eigenfunctions) of the GraField estimating equation based on the \texttt{G-matrix} under the proposed specialized nonparametric approximation scheme provides a systematic and unified framework for spectral graph analysis. As an application of this general formulation, we have shown how one can \emph{synthesize the well-known Laplacian, diffusion map, and modularity spectral algorithms} and view them as ``empirical'' spectral graph analysis methods. It is one of those rare occasions where one can witness the convergence of statistical, algorithmic and geometry-motivated computational models.
\end{rem}
\vspace{-.65em}

\section{Smoothed Spectral Graph Analysis} \label{sec:rspec}
Smoothness is a universal requirement for constructing credible nonparametric estimators. Spectral graph analysis is also no exception.
An improved smooth version of raw-empirical spectral graph techniques will be discussed, revealing a simple and straightforward statistical explanation of the \emph{origin} of regularized Laplacian techniques.
\subsection{High-dimensional Undersampled Regime}
Recall from Theorem \ref{thm:GEN} that the generalized matrix eigenvalue equation \eqref{eq:master} depends on the \emph{unknown} network and vertex probability mass functions. This leads us to the question of estimating the unknown distribution $P=(p_1,p_2,\ldots,p_n)$ (support size = size of the graph = $n$) based on $N$ sample, where $N=\sum_{i=1}^n d_i = 2|E|$. Previously (Theorems 5-7) we have used the \emph{unsmoothed} maximum likelihood estimate (MLE) $\tp(x;\cG)$ to construct our \emph{empirical} spectral approximation algorithms, which is the unique minimum variance unbiased estimator. Under standard asymptotic setup, where the dimension of the parameter space $n$ is fixed and the sample size $N$ tends to infinity, the law of large numbers ensures optimality of $\tp$. As a consequence, empirical spectral analysis techniques are expected to work quite well for dense graphs.
\vskip.55em

\emph{Estimation of probabilities from sparse data}. However, the raw $\tp$ is known to be strictly sub-optimal \citep{witten1991} and unreliable in the high-dimensional sparse-regime where $N/n = O(1)$ (i.e., when parameter dimension and the sample size are comparably large). This situation can easily arise for modern day \emph{large sparse graphs} where the ratio $N/n$ is small and there are many nodes with low degree, as is the case of degree sparsity. The naive MLE estimator can become unacceptably noisy (high variability) due to the huge size and sparse nature of the distribution. In order to reduce the fluctuations of ``spiky'' empirical estimates, some form of ``smoothing'' is necessary. The question remains: How to tackle this high-dimensional discrete probability estimation problem, as this directly impacts the quality of our nonparametric spectral approximation.

\vskip.55em

My main purpose in the next section is to describe one such promising technique for smoothing raw-empirical probability estimates, which is flexible and in principle can be applied to any sparse data.
\subsection{Spectral Smoothing}  We seek a practical solution for circumventing this problem that lends itself to fast computation. The solution, that is both the simplest and remarkably serviceable, is the Laplace/Additive smoothing \citep{laplace1951} and its variants, which excel in sparse regimes \citep{fienberg1973,witten1991}. The MLE and Laplace estimates of the discrete distribution $p(j;\cG_n)$ are respectively given by
\vskip.25em

Raw-empirical MLE estimates:  ~~~~~~~~~~~~$\tp(j;\cG_n)\,=\,\dfrac{d_j}{N}\,;$ ~~~~~~~~\\[.25em]
Smooth Laplace estimates:~~~~~~~~~~~~~~~~~~$\hp_{\tau}(j;\cG_n)\,=\,\dfrac{d_j + \tau}{N+n\tau}$~~~~~~~~~$(j=1,\ldots,n)$.
\vskip.25em

Note that the  smoothed distribution $\hp_{\tau}$ can be expressed as a convex combination of the empirical distribution $\tp$ and the discrete uninform distribution $1/n$
\beq \label{eq:lap} \hp_{\tau}(j;\cG_n)\,=\, \dfrac{N}{N+n\tau}\,\, \tp(j;\cG_n)\,+\,  \dfrac{n\tau}{N+n\tau} \, \left(\dfrac{1}{n}\right),~~~\eeq
which provides a Stein-type shrinkage estimator of the unknown probability mass function $p$. The shrinkage significantly reduces the variance, at the expense of slightly increasing the bias.
\vskip.45em

{\bf Choice of $\tau$}. The next issue is how to select the ``flattening constant'' $\tau$. The following choices of $\tau$ are most popular in the literature\footnote{For more details on selection of $\tau$ see \cite{fienberg1973} and references therein.}: \nocite{perks1947}
\[ \tau = \left\{ \begin{array}{ll}
         1 & \mbox{Laplace estimator};\\
        1/2 & \mbox{Krichevsky–Trofimov estimator}; \\
        1/n & \mbox{Perks estimator};\\
       \sqrt{N}/n &\mbox{Minimax estimator (under $L^2$ loss)}.\end{array} \right. \]
\vskip.25em

\begin{note}
Under increasing-dimension asymptotics, this class of estimator is often difficult to improve without imposing \emph{additional} smoothness constraints on
the vertex probabilities; see Bishop et al. (2007, Chapter 12)\nocite{bishop2007book}. The latter may not be a valid assumption as nodes of a graph offer \emph{no natural order} in general.
\end{note}

With this understanding, smooth generalizations of empirical spectral graph techniques will be discussed, which have a close connection with recently proposed spectral regularized techniques.
\subsection{Type-I Regularized Graph Laplacian}
Construct $\tau$-regularized smoothed empirical $\eta_{j;\tau}$ basis function by replacing the amplitude $p^{-1/2}(j)$ by $\hat{p}_\tau^{-1/2}(j)$ following \eqref{eq:lap}. Incorporating this regularized trial basis, we have the following modified \texttt{G-matrix} based linear algebraic estimating equation \eqref{eq:master}:
\beq \label{eq:m2}
\sum_j \te_{jk} \left[ \dfrac{\widetilde{p}(j,k)}{\sqrt{\hat p_\tau(j) \hat p_\tau(k)}}\,-\,\sqrt{\hat p_\tau(j)}\sqrt{\hat p_\tau(k)}\,-\,\la_k\delta_{jk} \right] ~=~0.
\eeq
\vskip.2em
\begin{thm}
The $\tau$-regularized block-pulse series based spectral approximation scheme is equivalent to representing or embedding discrete graphs in the continuous eigenspace of
\beq
\text{Type-I Regularized Laplacian} ~~=~~D_\tau^{-1/2}\,A\,D_\tau^{-1/2},~~~~~~~~~~~~
\eeq
where $D_\tau$ is a diagonal matrix with $i$-th entry $d_i+\tau$.
\end{thm}
\vskip.1em
\begin{note}
It is interesting to note that this exact regularized Laplacian formula was proposed by \cite{chaudhuri2012} and \cite{qin2013}, albeit from a very different motivation.
\end{note}
\subsection{Type-II Regularized Graph Laplacian}

\begin{thm}
Estimate the joint probability $p(j,k;\cG)$ by extending the formula given in \eqref{eq:lap} for the two-dimensional case as follows:
\beq \label{eq:2dlap}
\hp_{\tau}(j,k;\cG)\,\,=\,\, \dfrac{N}{N+n\tau}\,\, \tp(j,k;\cG) \,+\,  \dfrac{n\tau}{N+n\tau} \, \left(\dfrac{1}{n^2}\right),~~~~~~
\eeq
which is equivalent to replacing the original adjacency matrix by $A_\tau=A+(\tau/n) {\bf 1}{\bf 1}^{T}$. This modification via smoothing in the estimating equation \eqref{eq:m2} leads to the following spectral graph matrix
\beq \label{eq:rlp2}
\text{Type-II Regularized Laplacian} ~~=~~D_\tau^{-1/2}\,A_\tau\,D_\tau^{-1/2}.~~~~~~~~~~~~
\eeq
\end{thm}
\vskip.4em
\begin{note}
Exactly the same form of regularization of Laplacian graph matrix \eqref{eq:rlp2} was proposed by \cite{amini2013} as a fine-tuned empirical solution.
\end{note}
\begin{sig}
Recently there has been a lot of discussion on how to choose the spectral regularization parameter $\tau$. Using Davis–Kahan's theorem, \cite*{joseph2013} have proposed a data-driven technique called \emph{DKest}. Their method involves repeated eigen-decomposition of regularized Laplacian matrix (which is a \emph{dense} matrix) over a grid of values of $\tau$. There are two major concerns for \emph{DKest}: heavy computational burden and theoretical validity (which only holds for stochastic block models and its extensions). In contrast,  our analysis provides a strikingly simple recommendation without \emph{any additional} computational overhead (by connecting it to large sparse distribution smoothing), and thus is much easier to construct, implement, and use in practice. Similar to \cite{joseph2013}, the critical parameter of our theory is the ratio $N/n$--the average degree of the nodes, contrary to the previous results \citep{chaudhuri2012,qin2013} that relied on the minimum degree.
\end{sig}
\vspace{-.4em}
\subsection{Google's PageRank Method}
Smoothing of network and vertex probability distributions appearing in generalized matrix equation \eqref{eq:master} resulted in Type-I and Type-II regularized Laplacian methods. The third possibility is to directly smooth the conditional or transitional probability matrix to develop a regularized version of random walk Laplacian (which we call Type-III regularized Laplacian) method for large sparse graphs.
\vskip.55em

\emph{Smoothing conditional probability function}.  Consider a random walk on $\cG$ with transition probability $\calT(i,j;\cG)=\Pr(X_{t+1}=j|X_t=i) \ge 0$. Note that the smoothing \eqref{eq:2dlap} can equivalently be represented as
\beq \label{eq:Google1}
\calT_\tau(i,j;\cG)\,\,=\,\,\dfrac{A(i,j;\cG) + \tau/n}{d_i+\tau}\,\,=\,\, (1-\al_\tau) \calT(i,j;\cG)\,+\,\al_\tau \left(\dfrac{1}{n}\right),
\eeq
where the degree-adaptive regularization parameter $\al_\tau=\tau(d_i+\tau)^{-1}$. One can construct an equivalent and more simplified (non-adaptive) estimator by directly smoothing \textit{each row} of $\calT$ (as it is a row-stochastic or row-Markov matrix) via Stein-like shrinkage
\beq \label{eq:Google2} \calT_\al~=~(1-\al) \calT\,+\,\al\begin{bmatrix}
1/n & 1/n & \dots  & 1/n\\
\vdots & \vdots & \ddots & \vdots\\
1/n & 1/n & \dots  & 1/n
\end{bmatrix}, ~~~~0 < \al < 1. \eeq
Spectral analysis can now proceed on this regularized transition probability matrix (either \eqref{eq:Google1} or \eqref{eq:Google2}) by substituting empirical transition matrix $\widetilde{\calT}=D^{-1}A$.

\begin{sig}
The non-adaptive regularized random-walk Laplacian (transition matrix) $\calT_\al$ was introduced by Larry Page and Sergey Brin in 1996 \citep{GOOGLE} and is called the Google's PageRank matrix. What seems natural enough from our nonparametric smoothing perspective, is in fact known to be highly surprising \emph{adjustment}--the ``Teleportation trick''.
\vskip.05em
`\textit{Teleporting is the essential distinguishing feature of the PageRank random walk that had not appeared in the literature before}' -- \cite{vigna2005}.
\end{sig}

\begin{note} 
The Google's PageRank matrix (which is different from the empirical random walk Laplacian in an important way) is probably the most famous, earliest, and spectacular example of spectral regularization that was originally introduced to counter the problem of \textit{dangling nodes} (nodes with no outgoing edges)  by ``connecting'' the graph\footnote{From a statistical estimation viewpoint, this can be thought of as a way to escape from the ``zero-frequency problem'' for discrete probability estimation.}.
\end{note}
\vspace{-.25em}
Web is a huge heterogeneous sparse graph where dangling nodes (dead-ends) and disconnected components are quite common. As a result random walk can get stuck, and may cycle around an isolated set of pages. Smoothing allows the random walk to teleport to a web page uniformly at random (or adaptively based on degree-weights) whenever it hits a dead end. The steady state vector describes the long term visit rate--the PageRank score, computed via eigen-decomposition of $\calT_\al$.

\subsection{Other Generalizations}
The beauty of our statistical argument is that it immediately opens up several possibilities to construct \emph{new} types of spectral regularization schemes, which are otherwise hard to guess using previous understanding. Two such promising techniques are discussed here. The first one deals with
Stein smoothing with data-driven shrinkage parameter.

\begin{thm}
Under the sparse asymptotic setup where the domain size $n \rightarrow \infty$ at the same rate as the sample size $N$, the risk (expected squared-error loss) of smooth add-$\widehat{\tau}$ probability estimate $\widehat{p}_{\widehat{\tau}}$ \eqref{eq:lap} with the following data-driven choice
\beq \widehat{\tau}~=~\dfrac{N^2 - \sum_{i=1}^n d_i^2}{n\sum_{i=1}^n d_i^2 - N^2}\eeq
is uniformly smaller than the risk of the unsmoothed MLE estimator $\widetilde{p}=N^{-1}(d_1,\ldots,d_n)$. In addition, the risk of $\widehat{p}_{\widehat{\tau}}$ is uniformly smaller than the risk of the estimator formed by choosing $\tau=1/2$.
\end{thm}
For proof one can use the technique described in \cite{fienberg1973}. As a consequence of this result, we can further improve the spectral regularization algorithm by selecting $\tau$ in a data-driven way, with no extra computation.

\vskip.65em

While the Laplace or additive smoothing performs well in general, there are situations where they perform poorly \citep{gale1994}.
The Good-Turing estimator \citep{good1953} is often the next best choice which is given by
\beq \label{eq:GT} \widehat{p}_{{\rm GT}}(i;\cG)~=~\dfrac{\varpi_{d_i+1}}{\varpi_{d_i}} \cdot \dfrac{d_i+1}{N},  \eeq
where $\varpi_k$ denotes the number of nodes with degree $k$. An excellent discussion on this topic can be found in \cite{orlitsky2003}. One can plug in the estimate \eqref{eq:GT} into the equation \eqref{eq:m2} to generate new spectral graph regularization technique.
\begin{sig}
We wish to emphasize that our novelty lies in addressing the open problem of obtaining a \emph{rigorous interpretation and extension} of spectral regularization. To the best of our knowledge this is the \emph{first work} that provides a more formal and intuitive understanding of the \emph{origin} of spectral regularization. We have shown how the regularization naturally arises as a consequence of high-dimensional discrete data smoothing. In addition, this point of view allows us to select appropriate regularization parameter $\tau$ with \emph{no} additional computation.
\end{sig}

\begin{figure*}[!tth]
\centering
\vspace{2em}
\includegraphics[height=\textheight,width=.25\textwidth,keepaspectratio,trim=2.5cm 0cm 1.5cm 1.5cm]{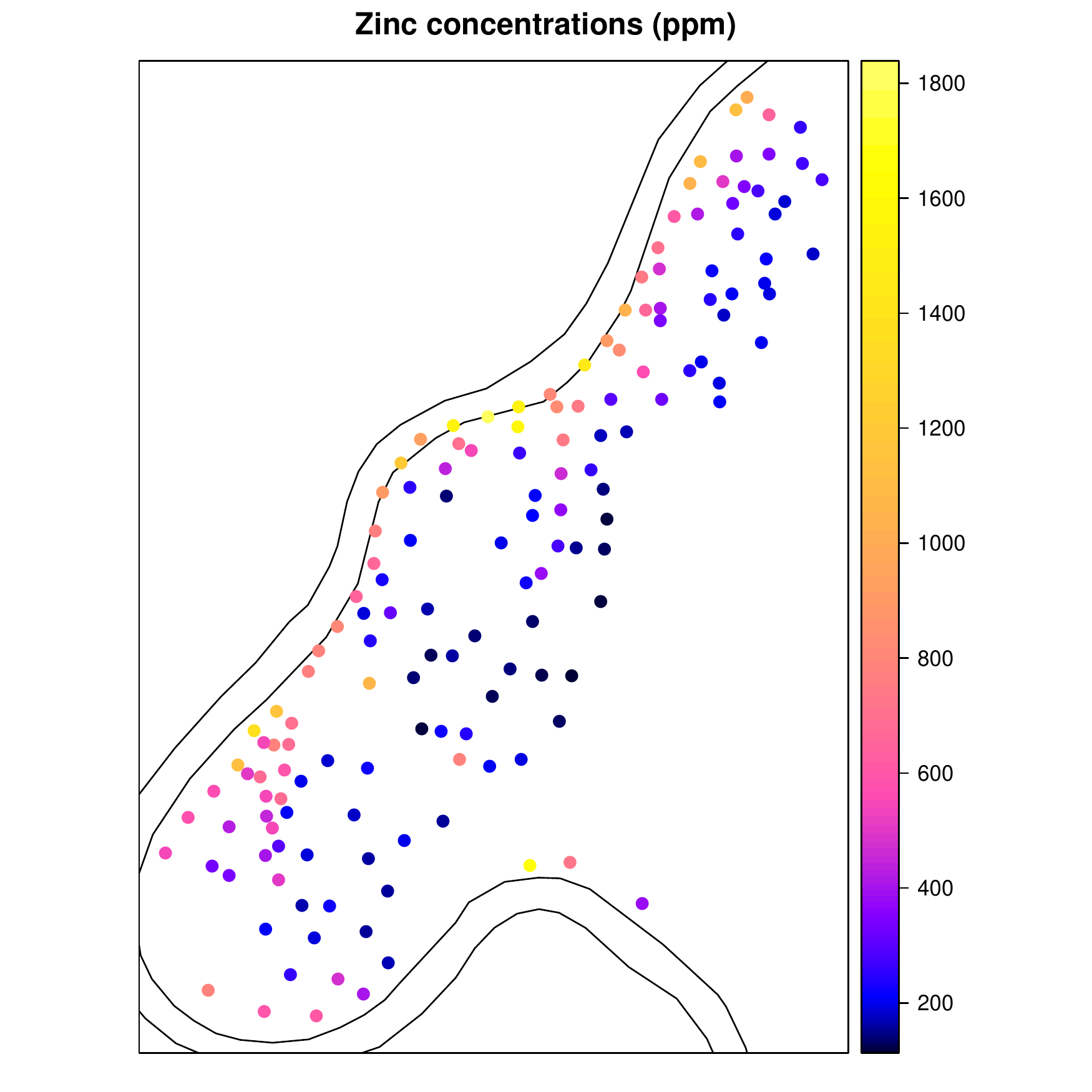}~
\includegraphics[height=\textheight,width=.25\textwidth,keepaspectratio,trim=2.5cm 0cm 1.5cm 1.5cm]{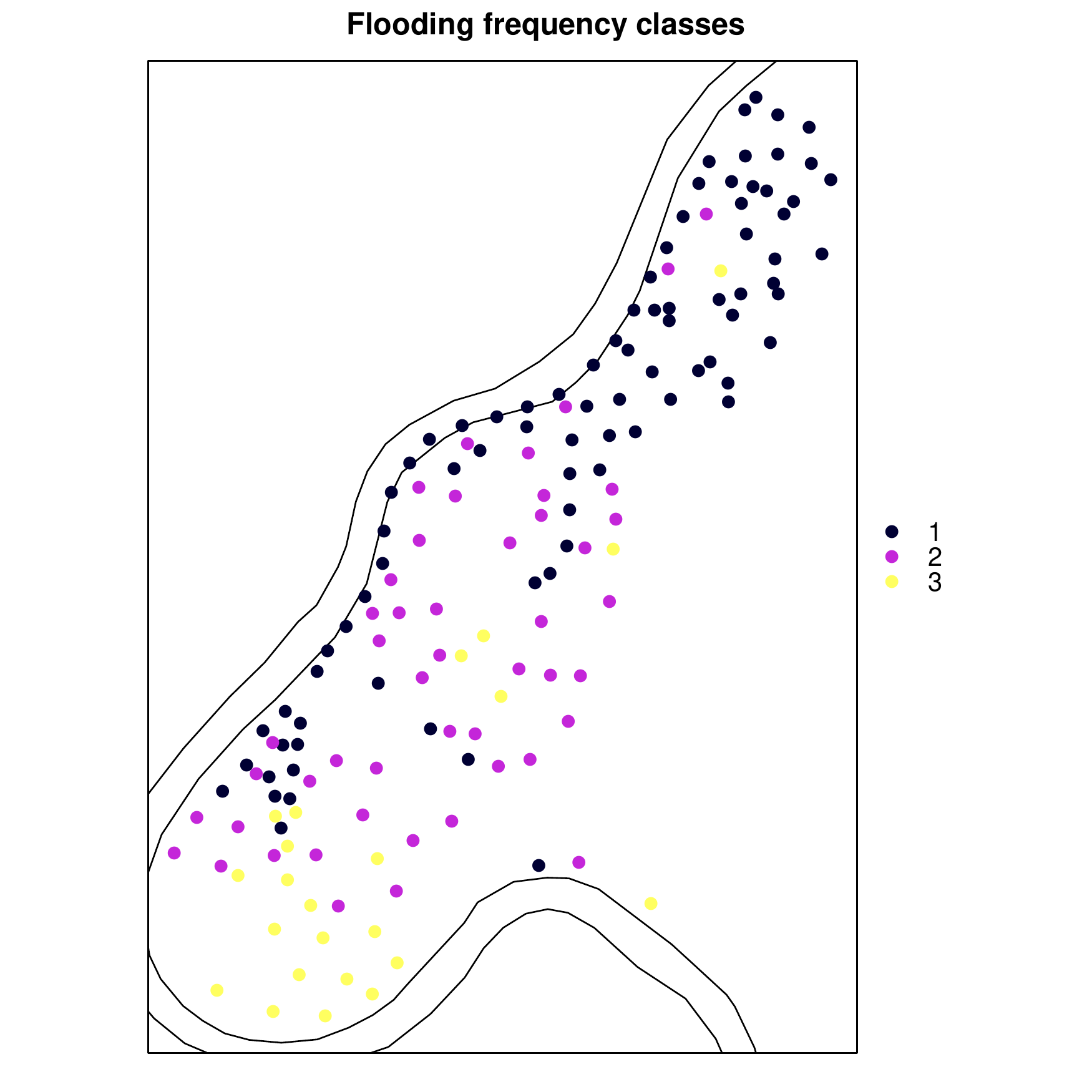}~
\includegraphics[height=\textheight,width=.25\textwidth,keepaspectratio,trim=2.5cm 0cm 1.5cm 1.5cm]{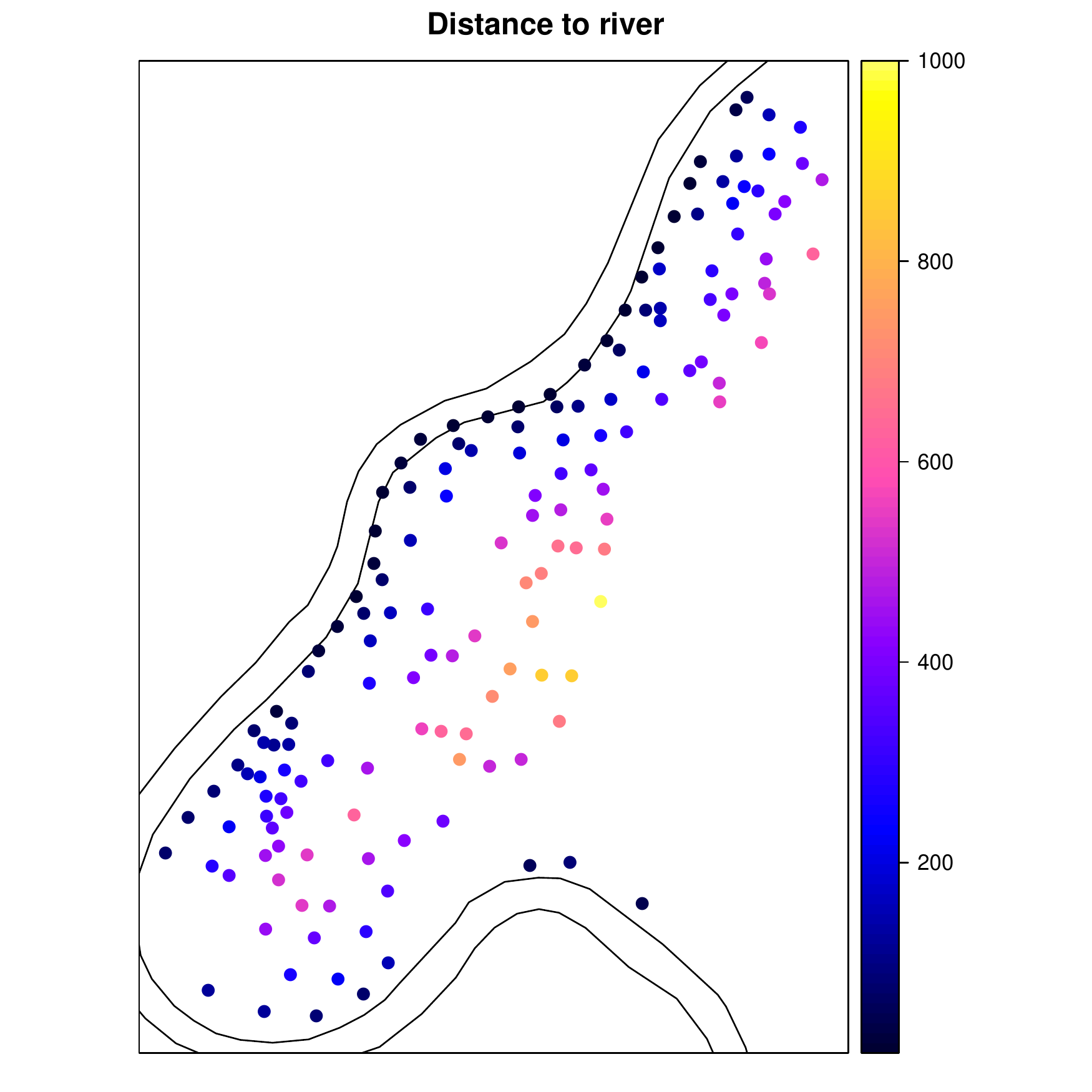}~
\includegraphics[height=\textheight,width=.25\textwidth,keepaspectratio,trim=2.5cm 0cm 1.5cm 1.5cm]{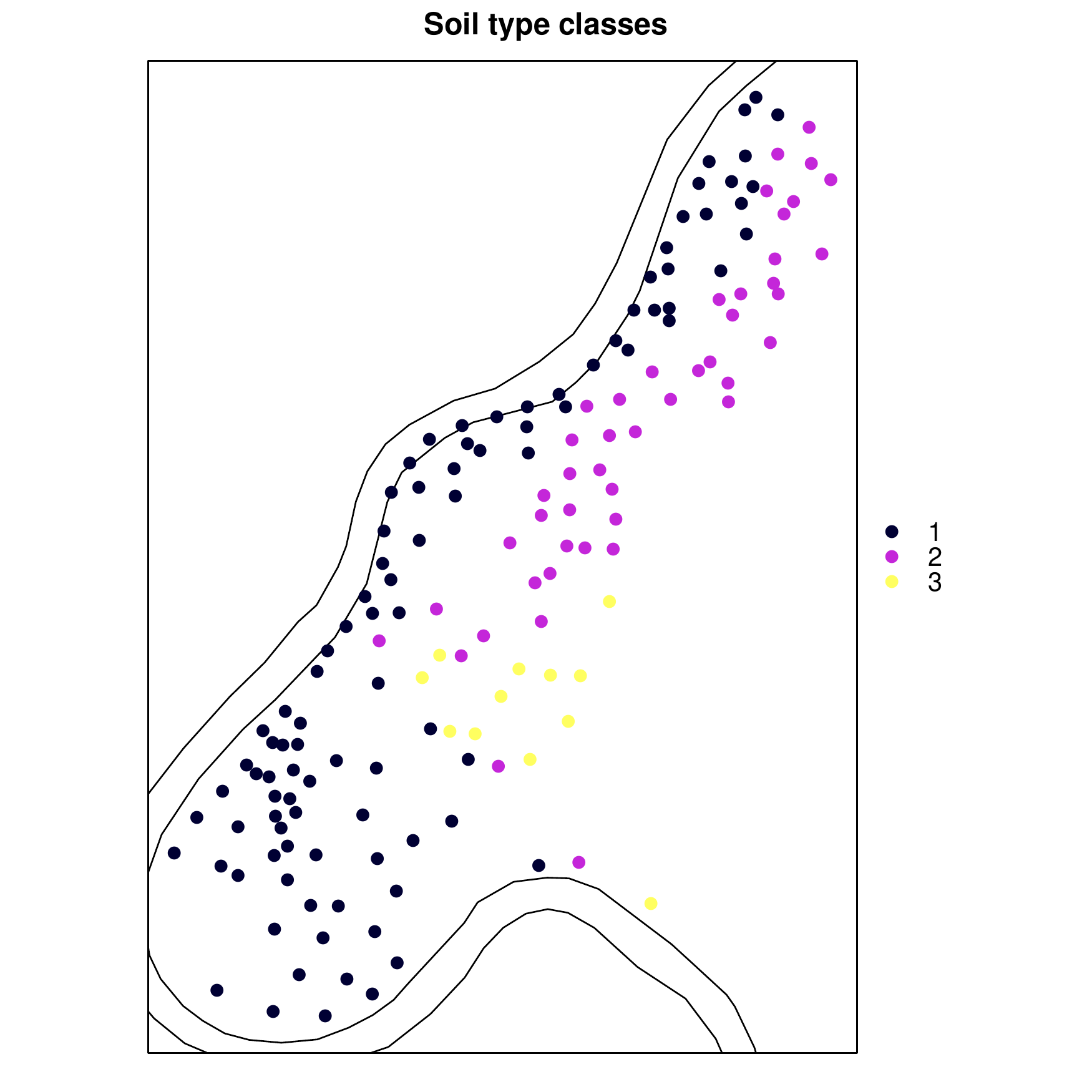}
\caption{The Meuse dataset consists of $n=155$ observations taken on a support of $15\times 15$ m from the top $0-20$ cm of alluvial soils in a $5\times 2$ km part of the right bank of the floodplain of the  the river Meuse, near Stein in Limburg Province (NL). The dependent variable $Y$ is the zinc concentration in the soil (in mg kg$^{-1}$), shown in the leftmost figure. The other three variables, flooding frequency class ($1$ = once in two years; $2$ = once in ten years; $3$ = one in $50$ years), distance to river Meuse (in metres), and soil type ($1$= light sandy clay; $2$ =  heavy sandy clay;  $3$ =silty light clay), are explanatory variables. \vspace{.4em}} \label{fig:meuse}
\end{figure*}

\section{Application to Graph Regression}
We study the problem of graph regression as an interesting application of the proposed nonparametric spectral analysis algorithm. Unlike traditional regression settings, here, one is given $n$ observations of the response and predictor variables over the graph. The goal is to estimate the regression function by properly taking into account the underlying graph-structured information along with the set of covariates.
\vskip.55em

\vskip.2em

\begin{figure*}[!tth]
\centering
\includegraphics[height=\textheight,width=.34\textwidth,keepaspectratio,trim=4cm 4cm 2cm 1.7cm]{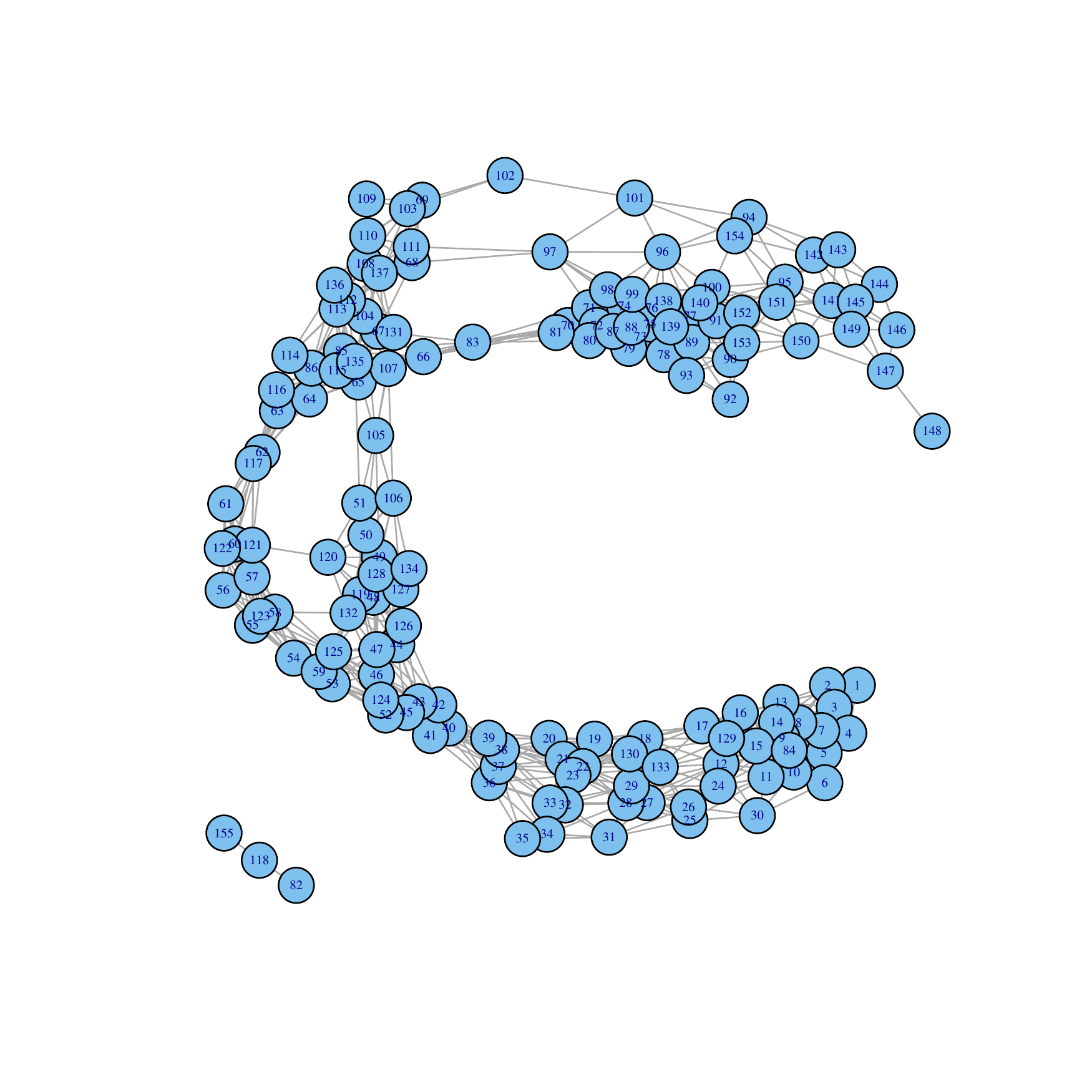}~~
\includegraphics[height=\textheight,width=.32\textwidth,keepaspectratio,trim=2cm 1.5cm 1cm 2cm]{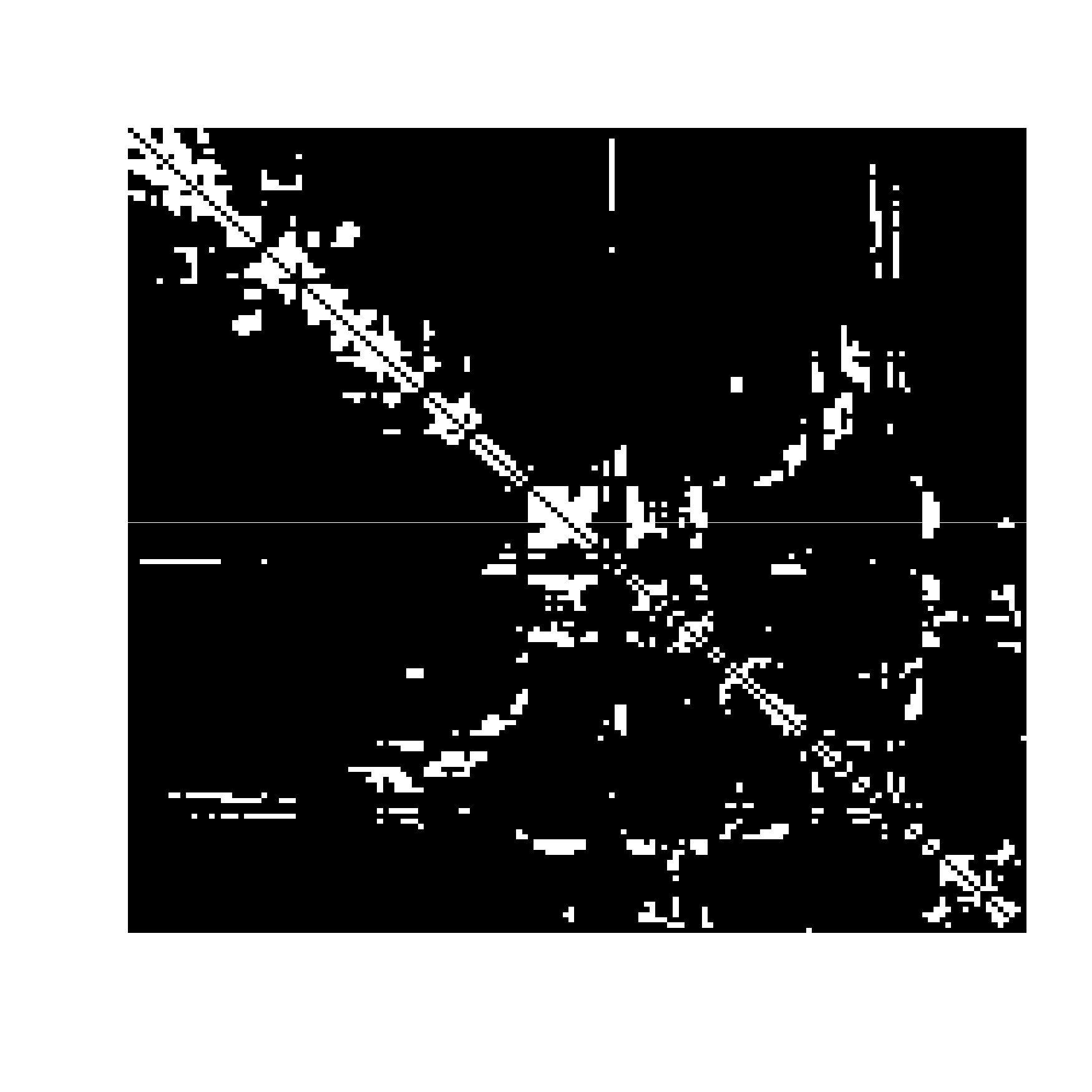}~~
\includegraphics[height=\textheight,width=.3\textwidth,keepaspectratio,trim=.5cm 1.5cm 2.5cm 2cm]{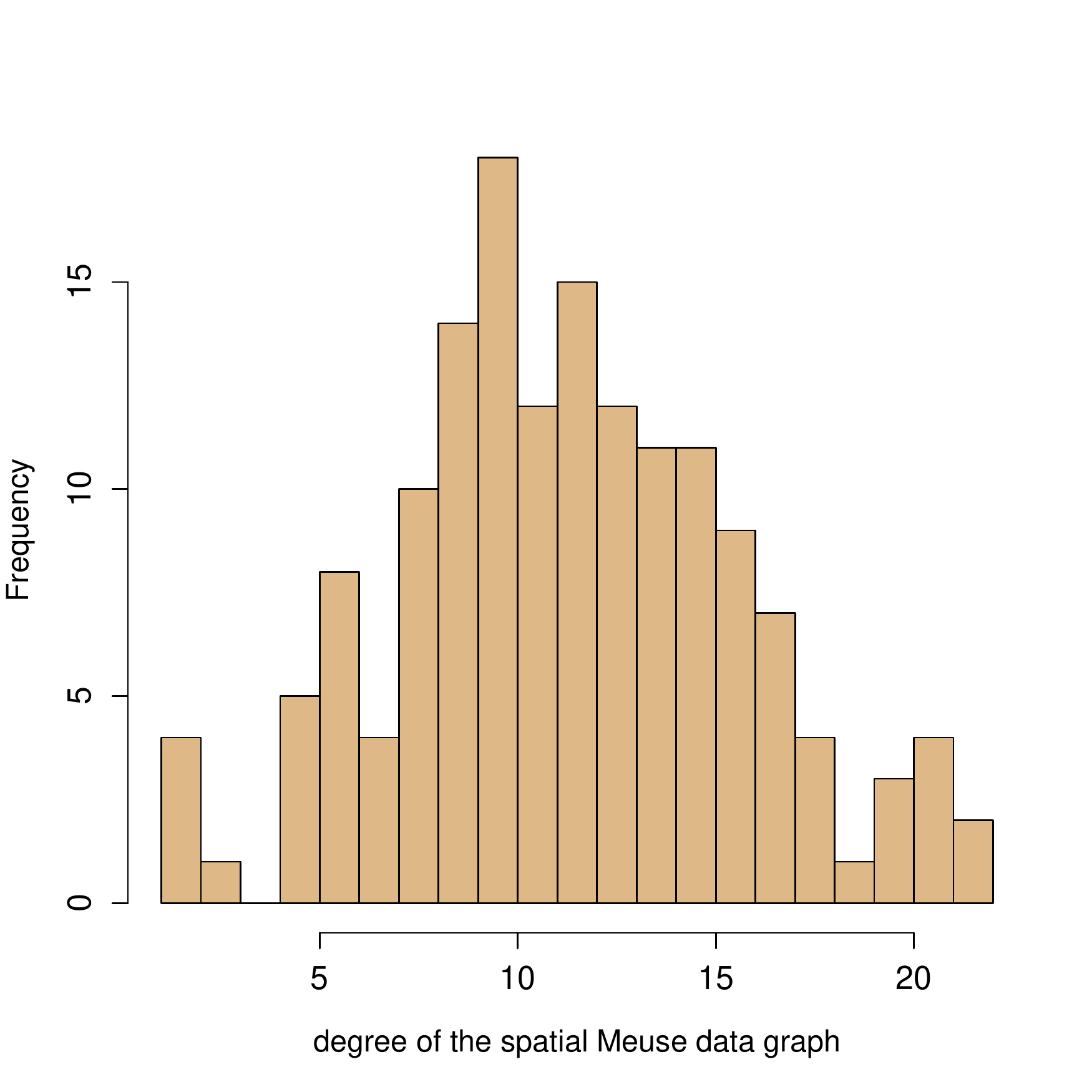} \vskip2em
\caption{From spatial to graph representation of the Meuse river data set. For each node we observe $\{y_i; x_{1i},x_{2i},x_{3i}\}_{i=1}^n$. We address the problem of approximating the regression function by simultaneously incorporating the effect of explanatory variables $X$ and the underlying spatial graph dependency.}\label{fig:meuseG}
\end{figure*}

{\bf Meuse Data Modeling}.   We apply our frequency domain graph analysis to address the spatial prediction problem. Fig \ref{fig:meuse} describes the Meuse data set, a well known geostatistical dataset. There is a considerable spatial pattern one can see from Fig \ref{fig:meuse}. We seek to estimate a smooth regression function of the dependent variable $Y$ (zinc concentration in the soil) via generalized spectral regression that can exploit this spatial dependency. The graph was formed according to the geographic distance between points based on the spatial locations of the observations. We convert spatial data into signal supported on the graph by connecting two vertices if the distance between two stations is smaller than a given coverage radius. The maximum of the first nearest neighbor distances is used as a coverage radius to ensure at least one neighbor for each node. Fig \ref{fig:meuseG}(C) shows the sizes of neighbours for each node ranging from $1$ to $22$. Three nodes (\# 82, 148, and 155) have only $1$ neighbour; additionally one can see a very weakly connected small cluster of three nodes, which is completely detached from the bulk of the other nodes. The reason behind this heterogeneous degree distribution (as shown in Fig \ref{fig:meuseG}) is the irregular spatial pattern of the Meuse data.

\vskip.5em
We model the relationship between $Y$ and spatial graph topology by incorporating nonparametrically learned spectral representation basis. Thus, we expand $Y$ in the eigenbasis of $\cd$ and the covariates for the purpose of smoothing, which effortlessly integrates the tools from harmonic analysis on graphs and conventional regression analysis. The model \eqref{eq:glasso} described in the following algorithm simultaneously promotes spatial smoothness and sparsity.
\begin{center}
Algorithm 2. \emph{Nonparametric Spectral Graph Regression}
\end{center}
\vspace{-.5em}
\medskip\hrule height .65pt
\vskip.55em
\texttt{Step 1.}  Input: We observe $\{y_i; x_{i1}, \ldots, x_{ip}\}$ at vertex $i$ of the graph $\cG=(V,A)$ with size $|V|=n$. The regularization parameters $\tau$.
\vskip.4em
\texttt{Step 2.}  Construct $\tau$-regularized block-pulse shape functions $\eta_{j;\tau} = \hat{p}_{j;\tau}^{-1/2} \ind(u_{j-1} < u \le u_{j})$ for $j=1,\ldots,n$.
\vskip.45em
\texttt{Step 3.}  Compute the smoothed \texttt{G-Matrix} with respect to trial bases $\{\eta_{k;\tau}\}_{1\le k \le n}$  \vspace{-.7em}
\[\mathcal{M}_{\tau}[j,k;\eta,\cG] ~= ~\Big\langle \eta_{j;\tau}, \int_0^1(\cd_n-1)\eta_{k;\tau} \Big\rangle_{L^2[0,1]}~~~ \text{for}~ j,k=1,\ldots,n.\]

\texttt{Step 4.} Construct smooth graph Fourier basis
\[\widehat{\phi}_{k;\tau}(u)=\sum_{j=1}^n u_{jk}\eta_{j;\tau},\,~ {\rm for}\,\, k=1,\ldots,n-1,\]
where $\mathcal{M}=U\Lambda U^{T} = \sum_k u_k \mu_k u_k^{T}$, and $u_{ij}$ are the elements of the singular vector $U=(u_1,\ldots,u_n)$. Construct spectral basis matrix $\Phi=\big[ \widehat \phi_{1;\tau}, \ldots, \widehat \phi_{k;\tau}\big] \in \cR^{n \times k}$ for the graph $\cG$.
\vskip.5em

\texttt{Step 5.} Construct combined graph regressor matrix $X_{\cG}=\big[ \Phi ; X  \big]$, where $X=[X_1,\ldots X_p] \in \cR^{n \times p}$ is the matrix of predictor variables.
\vskip.5em

\texttt{Step 6.} Solve for $\beta_{\cG}=(\be_\Phi, \be_X)^{ T}$
\beq \label{eq:glasso}
\widehat \be_{\cG} \,=\,\argmin_{\beta_{\cG} \in \cR^{k+p}} ~\| y - X_{\cG}\be_{\cG}\|_2^2\,\, + \,\,\la \| \be_{\cG} \|_1.
\eeq
\medskip\hrule height .65pt
\vskip1.5em

Algorithm 2 extends traditional regression to data sets represented over a graph. The proposed frequency domain smoothing algorithm efficiently captures the spatial graph topology via the spectral coefficients $\widehat \be_\Phi \in \cR^k$ and can be interpreted as covariate-adjusted discrete graph Fourier transform of the response variable $Y$. The $\ell_1$ sparsity penalty automatically selects the coefficients with largest magnitudes thus provides compression.

The following table shows that incorporating the spatial correlation in the baseline unstructured regression model using spectral orthogonal basis functions (which are estimated from the spatial graph with $k=25$) boosts the model fitting from $62.78\%$ to $80.47\%$, which is an improvement of approximately $18\%$.
\vskip1em
\begin{tabular}{ccc ccc cc}
\hline \\[-1em]
$X$ &  $\widetilde{\Phi}^0 + X$ &  $\widetilde{\Phi}_{\tau=1}^{\rm{I}} + X$ &  $\widetilde{\Phi}_{\tau=.5}^{\rm{I}} + X$ &  $\widetilde{\Phi}_{\tau=.28}^{\rm{I}} + X$ &  $\widetilde{\Phi}_{\tau=1}^{\rm{II}} + X$ &  $\widetilde{\Phi}_{\tau=.5}^{\rm{II}} + X$ &  $\widetilde{\Phi}_{\tau=.28}^{\rm{II}} + X$\\
\midrule
$62.78$ & $74.80$ & $78.16$ & $80.47$ & $80.43$ & $80.37$ & $80.45$ & $80.43$\\
\hline
\end{tabular}

\vskip1em
Here $\widetilde{\Phi}^0, \widetilde{\Phi}_\tau^{\rm{I}}$ and $\widetilde{\Phi}_\tau^{\rm{II}}$ are the graph-adaptive piecewise constant orthogonal basis functions derived respectively from the ordinary, Type-I regularized, and  Type-II regularized Laplacian matrix. Using other forms of regularization the performance does not change much, and thus we do not show it here. Our spectral method can also be interpreted as a kernel smoother, where the spatial dependency is captured by the discrete \texttt{GraField} $\cd$. We finally conclude that extension from traditional regression to graph-structured spectral regression significantly improves the model accuracy.
\section{Concluding Remarks}
What is the ``common core'' behind all existing spectral graph techniques?
Despite half a century of research, it remains one of the most formidable open issues, if not the core problem in modern spectral graph theory. The need for a unified theory become particularly urgent and pressing in the era of ``algorithm deluge'' in order to systematize the practice of graph data analysis.
\vskip.4em
As an attempt to shed some light on this previously unasked question, I have introduced a framework for statistical reasoning that appears to  offer a complete and autonomous description (totally free from the classical combinatorial or algorithmic linear algebra-based languages) with the intention of leading the reader to a broad understanding of how different techniques interrelate. The prescribed approach brings a fresh perspective by appropriately transforming the spectral problem into a nonparametric approximation and smoothing problem to arrive at new algorithms, as well as both practical and theoretical results. It was a great surprise for me to be able to deduce existing techniques from some underlying
basic principles in a self-consistent way.  At the end, there is always the hope that a comprehensive understanding gained through unification will inspire new tools, to attack some of the important open problems in this area.

%


\vskip1em
\noindent\textsc{\textbf{Acknowledgments}}: The author benefited greatly from many fruitful discussions and exchanges of ideas he had at the John W. Tukey 100th Birthday Celebration meeting, Princeton University September 18, 2015. The author is particularly grateful to Professor Ronald Coifman for suggesting the potential connection between our formulation and Diffusion map (random walk on graph). We also thank Isaac Pesenson for helpful questions and discussions.
\vskip1em


\newpage
\pagenumbering{arabic}
\renewcommand{\thepage} {S--\arabic{page}}
\renewcommand{\baselinestretch}{1.4}
\setlength{\parskip}{1.4ex}
\begin{center}
{\large {\bf \underline{Supporting Online Material}  \\[1em]  Unified Statistical Theory of Spectral Graph Analysis}}
\\[.2in]
Subhadeep Mukhopadhyay\\
Department of Statistical Science,  Temple University\\ Philadelphia, Pennsylvania, 19122, U.S.A.\\
\texttt{deep@temple.edu}\\[.2in]
\end{center}
To further demonstrate the potential application of smooth spectral graph algorithms, in this supplementary section we discuss the community detection problem that seeks to divides nodes into into $k$ groups (clusters), with larger proportion of edges inside the group (homogeneous) and comparatively sparser connections between groups to understand the large-scale structure of network. Discovering community structure is of great importance in many fields such as LSI design, parallel computing, computer vision, social networks, and image segmentation.

\section*{Graph Clustering}
In mathematical terms, the goal is to recover the graph signals (class labels) $y: V \mapsto \{1,2,\ldots,k\}$ based on the connectivity pattern or the relationship between the nodes. Representation of graph in the spectral or frequency domain via the nonlinear mapping $\Phi:\cG(V,E) \mapsto \cR^{m}$ using discrete KL basis of density matrix $\cd$ as the co-ordinate is the most important learning step in the community detection. This automatically generates spectral features $\{\phi_{1i},\ldots,\phi_{mi}\}_{1\le i \le n}$ for each vertex that can be used for building the distance or similarity matrix to apply k-means or hierarchical clustering methods. In our examples, we will apply k-means algorithm in the spectral domain, which seeks to minimizing the within-cluster sum of squares. In practice, often the most stable spectral clustering algorithms determine $k$ by spectral gap: $k=\arg\!\max_j|\la_j - \la_{j+1}| + 1$.

\vspace{-.3em}
\begin{center}
Algorithm 3. \emph{Nonparametric Spectral Graph Partitioning}
\end{center}
\vspace{-.5em}
\medskip\hrule height .65pt
\vskip.55em
1. Input: The adjacency matrix $A \in \cR^{n \times n}$. Number of clusters $k$. The regularization parameter $\tau$.
\vskip.25em
2. Estimate the top $k-1$ spectral connectivity profile for each node $\{\widetilde{\phi}_{1i;\tau}, \ldots,   \widetilde{\phi}_{(k-1)i;\tau}\}$ using Algorithm 2. Store it in $\Phi \in \cR^{n \times k-1}$.
\vskip.25em
3. Apply k-means clustering by treating each row of $\Phi$ as a point in $\cR^{k-1}$.
\vskip.25em
4. Output: The cluster assignments of $n$ vertices of the graph $C_1,\ldots,C_k$.
\vskip.2em
\medskip\hrule height .65pt
\vskip1em

In addition to graph partitioning, the spectral ensemble $\{\la_k,\phi_k\}_{1\le k \le m}$ of $\cd$ contain a wealth of information on the graph structure. For example, the quantity $1-\widetilde{\la}_1(\cG;\eta)$ for the choice of $\{\eta_k\}$ to be normalized top hat basis \eqref{eq:basis}, is referred to as the algebraic connectivity, whose magnitude reflects how well connected the overall graph is. The kmeans clustering after spectral embedding $\{\widetilde{\phi}_j(\cG;\eta)\}_{1\le j \le k-1}$ finds approximate solution to the NP-hard combinatorial optimization problem based on the normalized cut \citep{shi2000} by relaxing the discreteness constraints into one that is continuous (sometimes known as  spectral relaxation).
\vskip.5em

{\bf Data and Results}.  We investigate four well-studied real-world networks for community structure detection based on $7$ variants of spectral clustering methods.

{\bf Example A} [Political Blog data, \cite{adamic2005}]  The data, which contains $1222$ nodes and $16,714$ edges, were collected over a period of two months preceding the U.S. Presidential Election of $2004$ to study how often the political blogs refer to one another.  The linking structure  of the political blogosphere was constructed by identifying whether a URL present on the page of one blog references another political blog (extracted from blogrolls). Each blog was manually labeled as liberal or conservative by \cite{adamic2005}, which we take as ground truth. The goal is to discover the community structure based on these blog citations, which will shed light on the polarization in political blogs.

Table \ref{tab:simulation} shows the result of applying the spectral graph clustering algorithm on this political web-blog data. The un-regularized Laplacian performs very poorly, whereas as both type-I/II regularized versions give significantly better results. The misclassification error drops from $47.95\%$ to $4.7\%$ because of regularization. To better understand why regularization plays a vital role, consider the degree distribution of the web-blog network as shown in the bottom panel of Figure 2. It clearly shows the presence of a large number of low-degree nodes, which necessitates the smoothing of high-dimensional discrete probability $p_1,\ldots,p_{1222}$. Thus, we perform the kmeans clustering after projecting the graph in the Euclidean space spanned by Laplace smooth KL spectral basis $\{\widetilde{\phi}_{k;\tau}\}$. Regularized spectral methods correctly identify two dense clusters: liberal and conservative blogs, which rarely links to a blog of a different political leaning, as shown in the middle panel of Fig \ref{fig:app}.

\begin{table}[t]
\caption{\textit{We report \% of misclassification error. We compare following seven different Laplacian variants. $K$ denotes the number of communities. }}
\vskip1em
\centering
\setlength{\tabcolsep}{5.5pt}
\def\arraystretch{1.128}%
\begin{tabular}{cccccccccc }
\toprule
   &&&\multicolumn{3}{c}{Type-I Reg. Laplacian} & &\multicolumn{3}{c}{Type-II Reg. Laplacian}\\ \cmidrule{4-6} \cmidrule{8-10}
Data & K & Laplacian &$\tau=1$&$\tau=1/2$&$\tau=\sqrt{N}/n$&  &$\tau=1$&$\tau=1/2$&$\tau=\sqrt{N}/n$\\ \midrule
PolBlogs & $2$ &47.95\%& 4.9\%  & 4.8\%  & 5.4\%  && 4.8\% & 4.7\%  &5.4\%\\

\midrule

Football & $11$  & 11.3\% & 7.83\%  & 6.96\%  & 6.96\% &  & 6.96\%  & 7.83\%  & 7.83\%\\

\midrule

MexicoPol &$2$ & 17.14\% & 14.2\% & 14.2\%& 14.2\% && 14.2\%& 14.2\%& 14.2\%\\

\midrule
Adjnoun &$2$ & 13.4\% & 12.5\%  & 12.5\% & 12.5\%& & 12.5\%& 12.5\%& 12.5\%\\
%


\bottomrule
\end{tabular}
\vskip1em
\label{tab:simulation}
\end{table}

{\bf Example B} [US College Football, \cite{NG}] The American football network (with $115$ vertices, $615$ edges) depicts
the schedule of football games between NCAA Division IA colleges during the regular season of Fall $2000$. Each node represents a college team (identified by their college names) in the division, and  two teams are linked if they have played each other that season.  The teams were divided into $11$ ``conferences'' containing around 8 to 12 teams each, which formed actual communities. The teams in the same conference played more often compared to the other conferences, as shown in the middle panel of Fig \ref{fig:app}. A team played on average 7 intra- and 4 inter-conference games in the season. Inter-conference play is not uniformly distributed; teams that are geographically close to one another but belong to different conferences are
more likely to play one another than teams separated by large geographic distances. As the communities are well defined, the American football network provides an excellent real-world benchmark for testing community detection algorithms.

Table \ref{tab:simulation} shows the performance of spectral community detection algorithms to identify the $11$ clusters in the American football network data. The regularization boosts the performance by $3$-$4\%$. In particular, $\tau=1/2$ and $\sqrt{N}/n$ produces the best result for Type-I regularized Laplacian, while $\tau=1$ exhibits the best performance for Type-II regularized Laplacian.
\vskip.5em

{\bf Example C} [The Political Network in Mexico, \cite{gil1996}] The data (with $35$ vertices and $117$ edges) represents the complex relationship between politicians in Mexico (including presidents and their close associates). The edge between two politicians indicates a significant tie, which can either be political, business, or friendship. A classification of the politicians according to their professional background (1 - military force, 2 - civilians: they fought  each other for power) is given. We use this information to compare our $7$ spectral community detection algorithms.

\begin{figure*}[thb]
\centering
\vspace{1em}
\includegraphics[height=\textheight,width=.37\textwidth,keepaspectratio,trim=4cm 3.5cm 5cm 4cm]{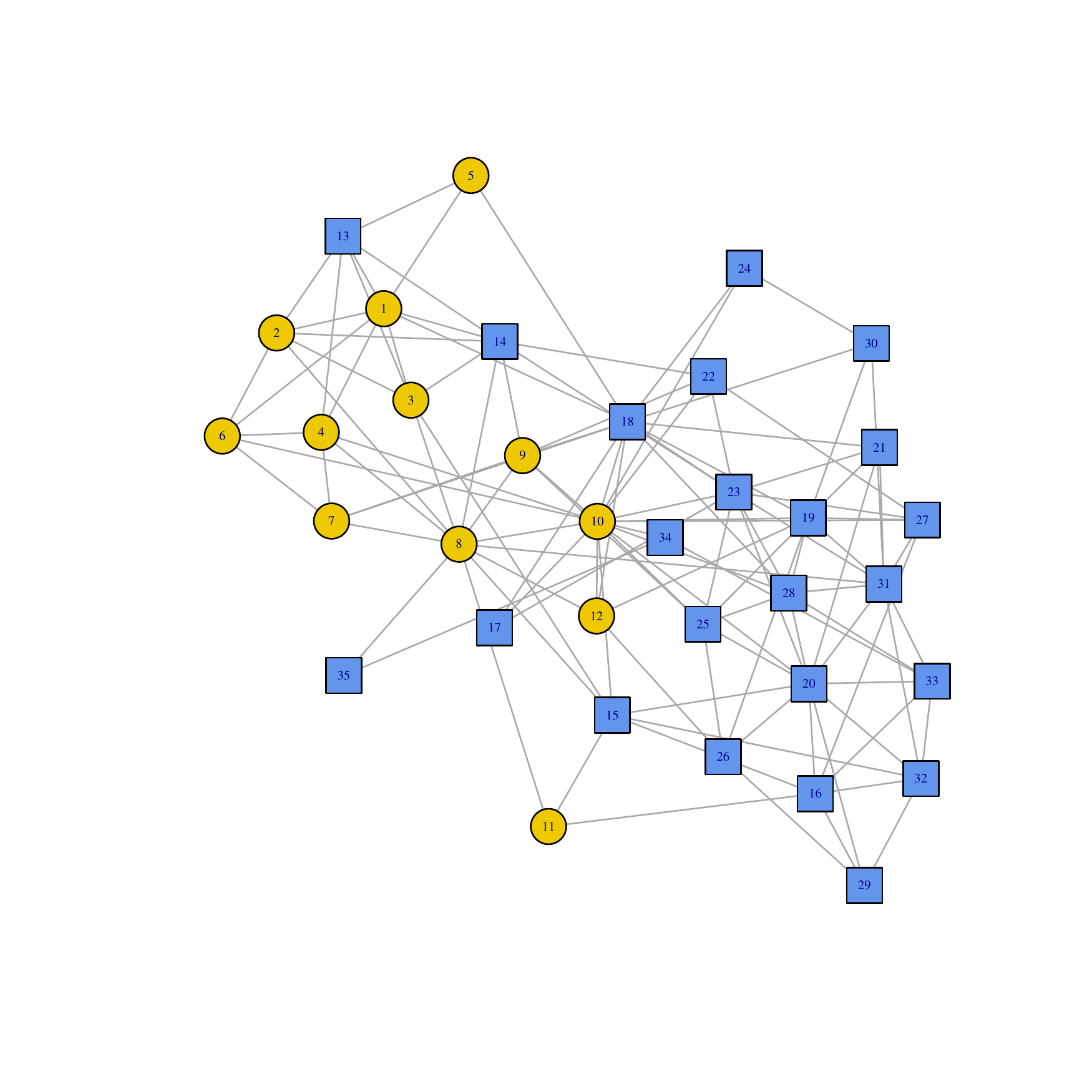}\\[1em]
\caption{ Mexican political network. Two different colors (golden and blue) denotes the two communities (military and civilians).} \label{fig:app2}
\end{figure*}

Although this is a ``small'' network, challenges arise from the fact that the two communities cannot be separated easily due to the presence of a substantial number of between-community edges, as depicted in Figs \ref{fig:app2} and \ref{fig:app}. The degree-sparsity is also evident from Fig \ref{fig:app} (bottom-panel). Table \ref{tab:simulation} compares seven spectral graph clustering methods. Regularization yields $3\%$ fewer  misclassified nodes. Both the type-I and II regularized Laplacian methods for all the choices of $\tau$ produce the same result.

\vskip1em
{\bf Example D} [Word Adjacencies, \cite{newman2006}] This is a adjacency network (with $112$ vertices and $425$ edges) of common adjectives and nouns in the novel \emph{David Copperfield} by English $19$th century writer Charles Dickens. The graph was constructed by \cite{newman2006}. Nodes represent the $60$ most commonly occurring adjectives and nouns, and an edge connects any two words that appear adjacent to one another at any point in the book. Eight of the words never appear adjacent to any of the others and are excluded from the network, leaving a total of $112$ vertices. The goal is to identify which words are adjectives and nouns from the given adjacency network.

\vskip.25em

Note that typically adjectives occur next to nouns in English. Although it is possible for adjectives to occur next to other adjectives
(e.g., ``united nonparametric statistics'') or for nouns to occur next to other nouns (e.g., ``machine learning''), these juxtapositions
are less common. As expected, Fig \ref{fig:app} (middle panel) shows an approximately bipartite connection pattern among the nouns and adjectives.

\vskip.25em
A degree-sparse (skewed) distribution is evident from the bottom right of Fig \ref{fig:app}. We apply the seven spectral methods and the result is shown in Table \ref{tab:simulation}. The traditional Laplacian yields a $13.4\%$ misclassification error for this dataset. We get better performance (although the margin is not that significant) after spectral regularization via Laplace smoothing.

\begin{figure*}[!thb]
\centering
\includegraphics[height=\textheight,width=.24\textwidth,keepaspectratio,trim=2cm 1.5cm 1cm 2cm]{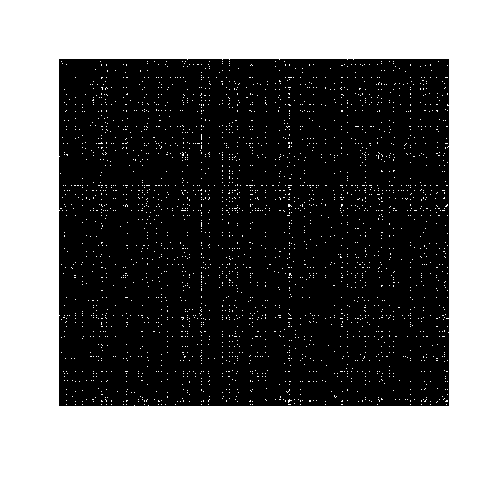}
\includegraphics[height=\textheight,width=.24\textwidth,keepaspectratio,trim=1.75cm 1.5cm 1cm 2.5cm]{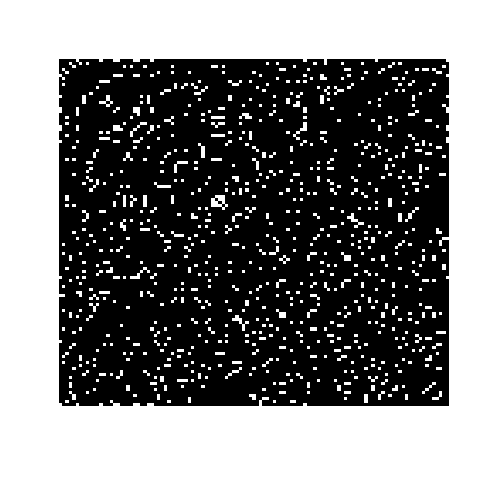}
\includegraphics[height=\textheight,width=.24\textwidth,keepaspectratio,trim=1.75cm 1.5cm 1cm 2.5cm]{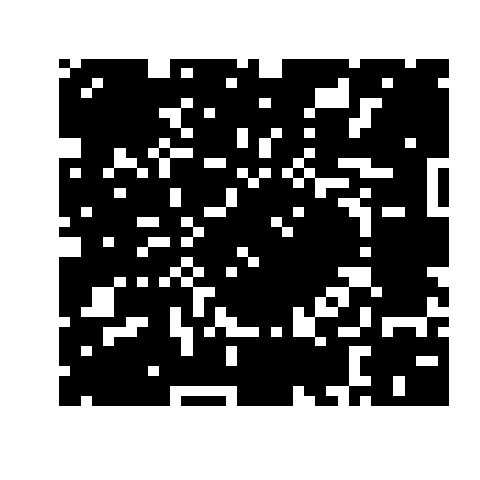}
\includegraphics[height=\textheight,width=.24\textwidth,keepaspectratio,trim=1cm 1.5cm 2cm 2cm]{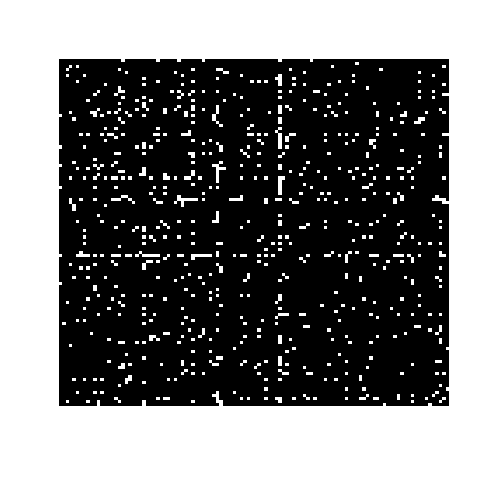}\\[2em]

\includegraphics[height=\textheight,width=.24\textwidth,keepaspectratio,trim=2cm 1.5cm 1cm 2cm]{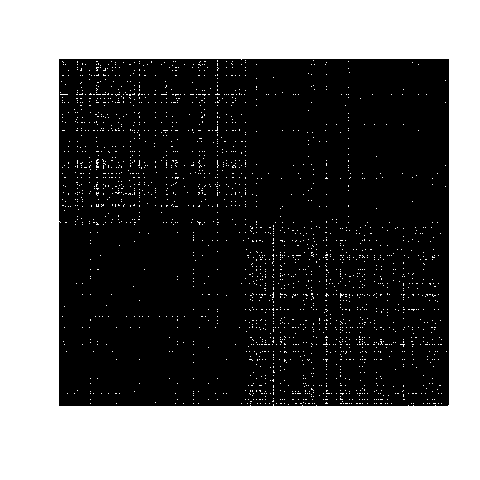}
\includegraphics[height=\textheight,width=.24\textwidth,keepaspectratio,trim=1.75cm 1.5cm 1cm 2.5cm]{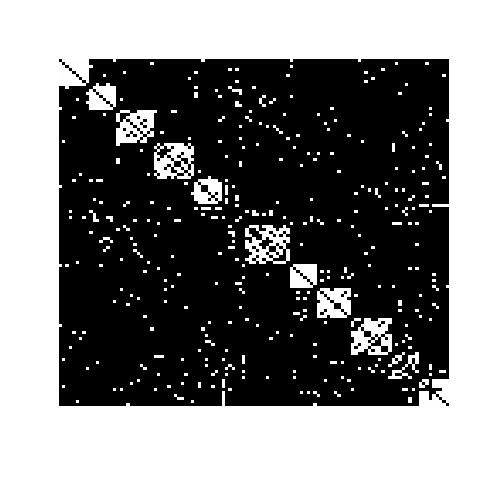}
\includegraphics[height=\textheight,width=.24\textwidth,keepaspectratio,trim=1.75cm 1.5cm 1cm 2.5cm]{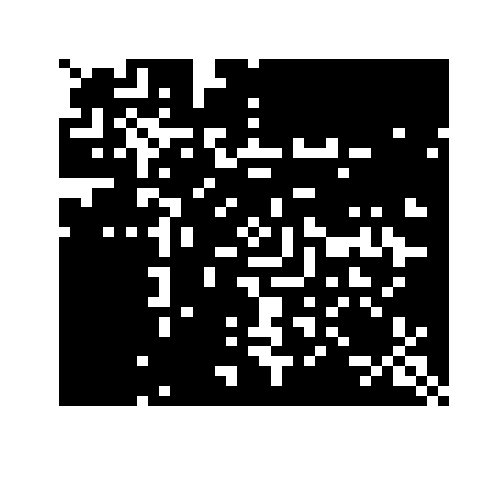}
\includegraphics[height=\textheight,width=.24\textwidth,keepaspectratio,trim=1cm 1.5cm 2cm 2cm]{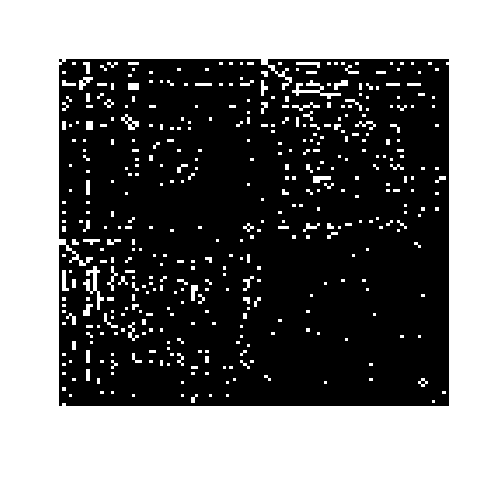}\\[2em]

\includegraphics[height=\textheight,width=.24\textwidth,keepaspectratio,trim=2cm 1.5cm 1cm 2cm]{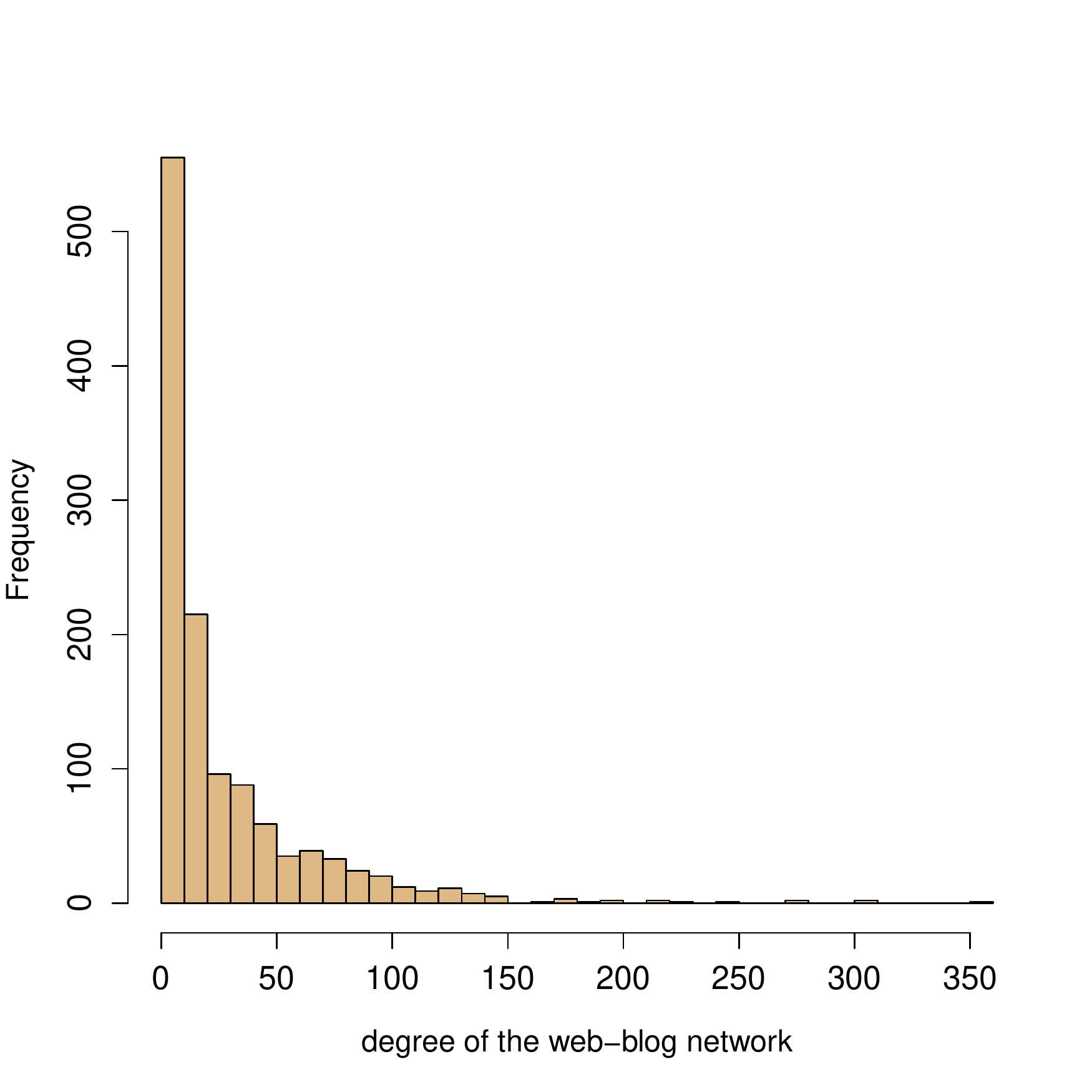}
\includegraphics[height=\textheight,width=.24\textwidth,keepaspectratio,trim=1.75cm 1.5cm 1cm 2.5cm]{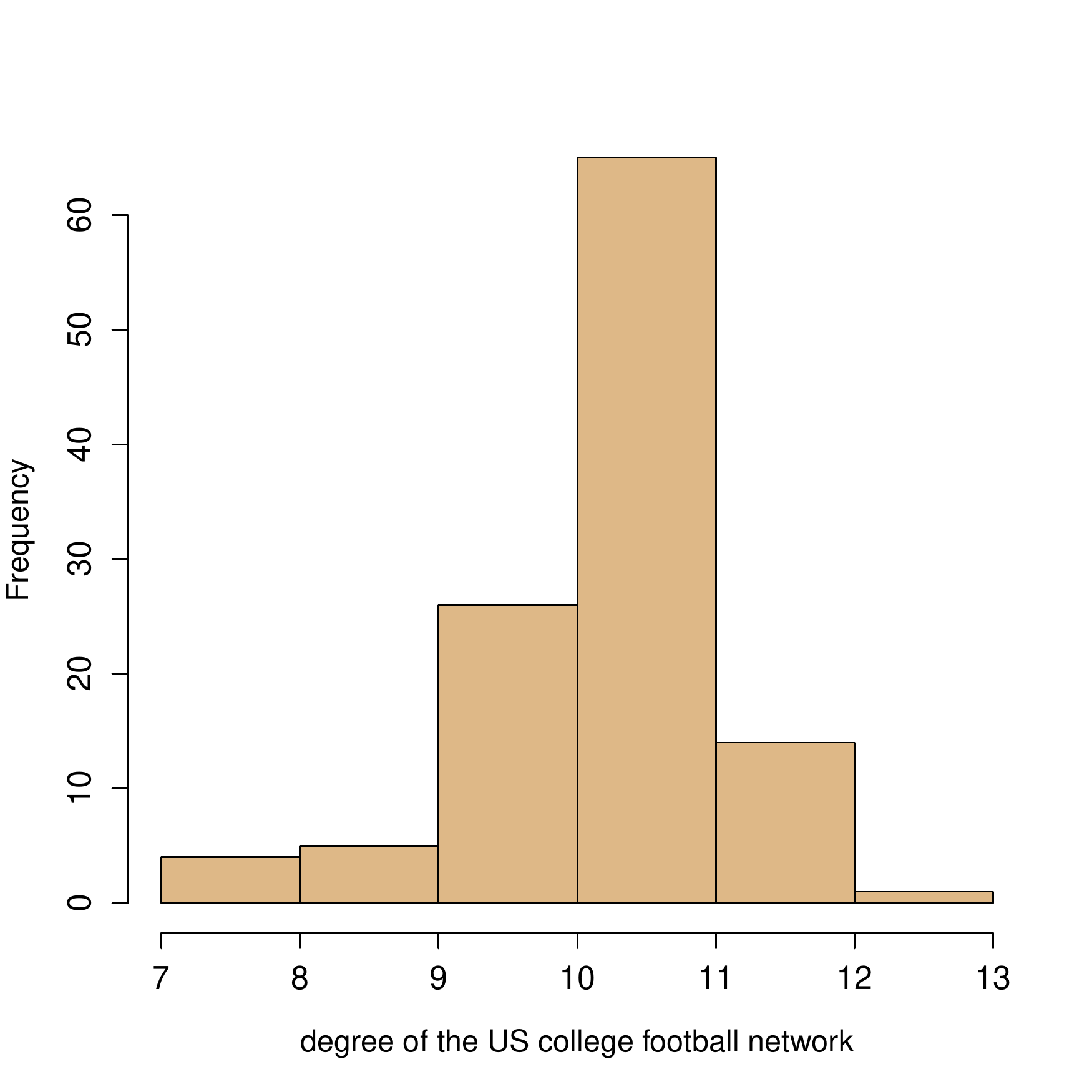}
\includegraphics[height=\textheight,width=.24\textwidth,keepaspectratio,trim=1.75cm 1.5cm 1cm 2.5cm]{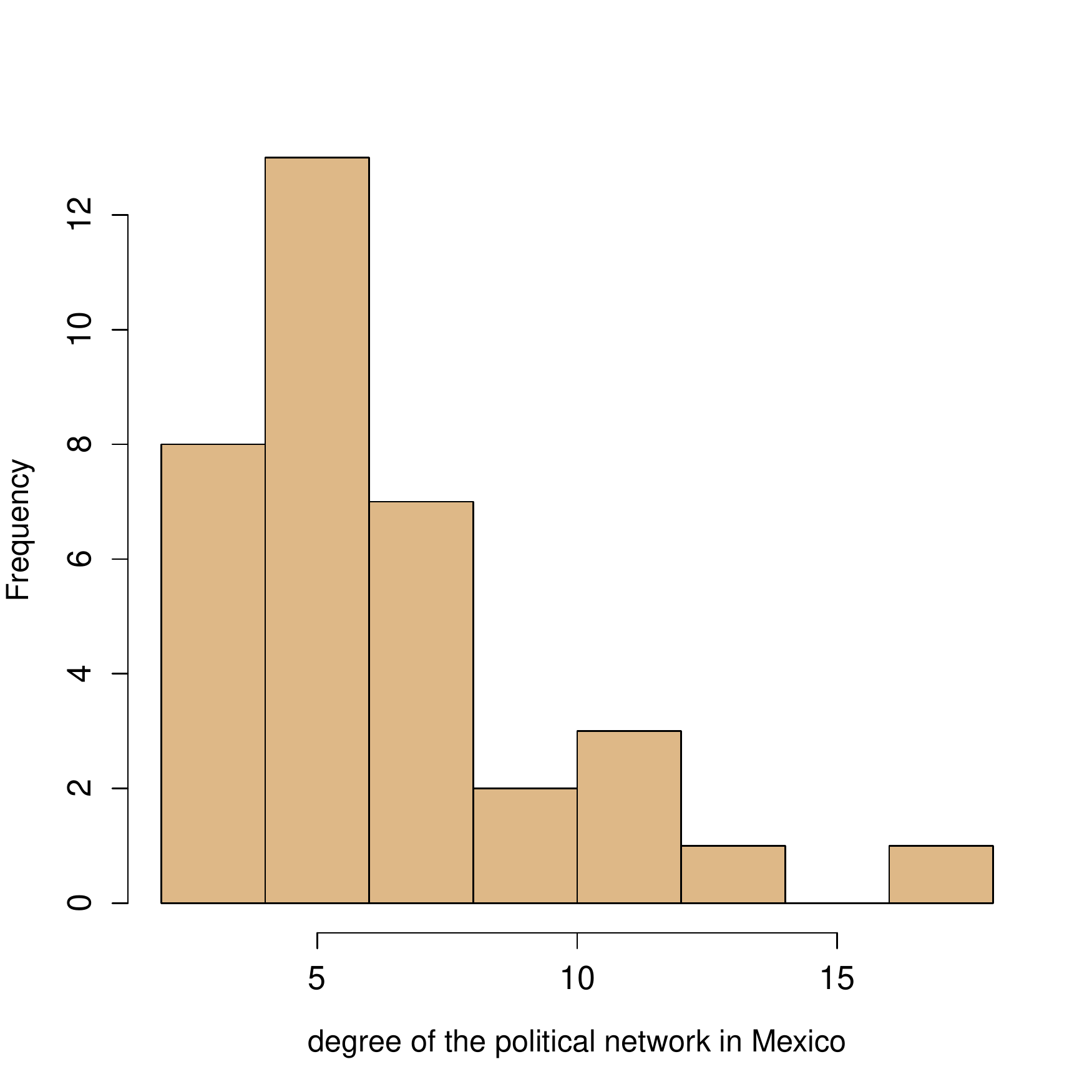}
\includegraphics[height=\textheight,width=.24\textwidth,keepaspectratio,trim=1cm 1.5cm 2cm 2cm]{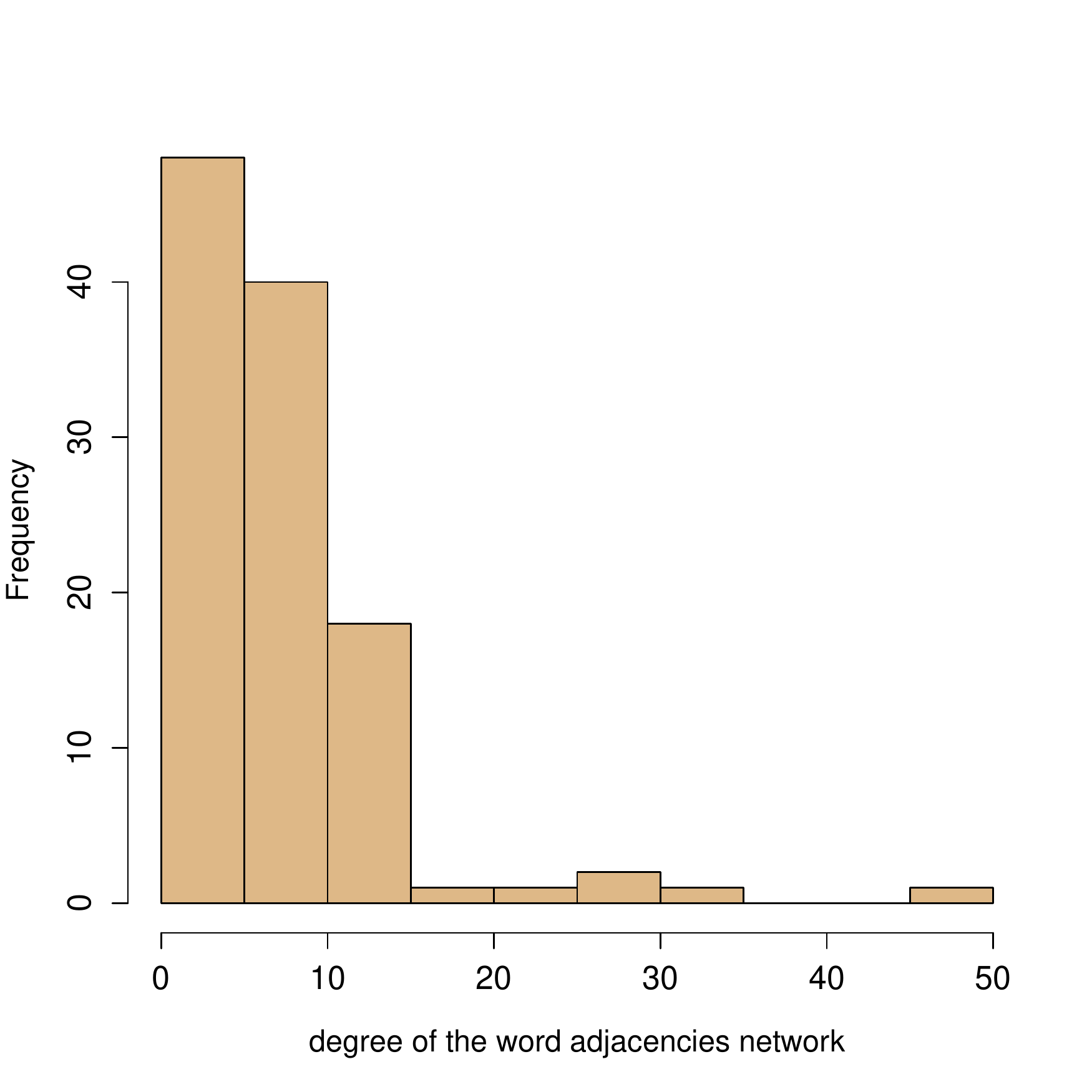}\\[2em]
\caption{The columns denote the 4 datasets corresponding to the political blog, US college football, politicians in Mexico, and word adjacencies networks. The first two rows display the un-ordered and ordered adjacency matrix, and the final row depicts the degree distributions.} \label{fig:app}
\vskip1em
\end{figure*}

\end{document}